\input ./style/arxiv-general.cfg
\documentclass[aos,MSNbibl,seceqn,dvips]{arximspdf}
\makeatletter
   \@ifpackageloaded{graphicx}{}{\usepackage{graphicx}}
\makeatother

\doi{10.1214/15-AOS1318}
\volume{43}
\issue{4}
\pubyear{2015}
\firstpage{1617}
\lastpage{1646}
\docsubty{FLA}

\makeatletter
\newcommand{\rrvert}{\vert}
\newcommand{\rrVert}{\Vert}
\newcommand{\llvert}{\vert}
\newcommand{\llVert}{\Vert}
\newcommand{\vvvert}{|\!|\!|}
\newcommand{\duvvvert}{\big|\!\big|\!\big|}
\newcommand{\cC}{{\mathcal C}}
\newcommand{\cD}{{\mathcal D}}
\newcommand{\cF}{{\mathcal F}}
\newcommand{\cH}{{\mathcal H}}
\newcommand{\cM}{{\mathcal M}}
\newcommand{\cN}{{\mathcal N}}
\newcommand{\cO}{{\mathcal O}}
\newcommand{\cS}{{\mathcal S}}
\newcommand{\cT}{{\mathcal T}}
\newcommand{\cX}{{\mathcal X}}
\newcommand{\cZ}{{\mathcal Z}}
\newcommand{\bE}{\mathbb E}
\newcommand{\bL}{{\mathbb L}}
\newcommand{\bN}{{\mathbb N}}
\newcommand{\bP}{{\mathbb P}}
\newcommand{\bR}{{\mathbb R}}
\newtheorem{theorem}{Theorem}
\newtheorem{lemma}{Lemma}
\newproclaim{defi}{Definition}
\newcommand{\argmin}{\mathop{\operatorname{arg}\min}}
\newcommand{\ind}{\mathrm{1}}
\newcommand{\R}{R}
\newcommand{\D}{G}
\newcommand{\Rn}{\widehat{R}}
\newcommand{\Dn}{\widehat{G}}
\newcommand{\f}{\theta}
\newcommand{\fn}{\widehat{\theta}}
\newcommand{\RC}{\mathcal{W}}
\newcommand{\DC}{\nabla\mathcal{W}}
\newcommand{\RCn}{\widehat{\mathcal{W}}}
\newcommand{\DCn}{\nabla\widehat{\mathcal{W}}}
\newcommand{\bc}{\mathbf{c}}
\newcommand{\bcn}{\widehat{\mathbf{c}}}
\newcommand{\majorant}{\mathcal{M}}
\newcommand{\bias}{B}
\newcommand{\estim}{\widehat{f}}
\makeatother

\begin{document}
\begin{frontmatter}

\title{Bandwidth selection in kernel empirical risk minimization via
the gradient}
\runtitle{Bandwidth selection via the gradient}

\begin{aug}
\author[A]{\fnms{Micha\"{e}l}~\snm{Chichignoud}\thanksref{T1}\ead[label=e1]{chichignoud@stat.math.ethz.ch}}
\and
\author[B]{\fnms{S\'{e}bastien}~\snm{Loustau}\corref{}\ead[label=e2]{loustau@math.univ-angers.fr}}
\runauthor{M. Chichignoud and S. Loustau}
\affiliation{ETH Z\"{u}rich and University of Angers}
\address[A]{Seminar Fuer Statistics\\
ETH Z\"urich\\
R\"amistrasse 101\\
CH-8092 Z\"urich\\
Switzerland\\
\printead{e1}}
\address[B]{LAREMA\\
Universit\'e d'Angers\\
2 Bvd Lavoisier\\
49045 Angers Cedex\\
France\\
\printead{e2}}
\end{aug}
\thankstext{T1}{Supported in part as member of the German--Swiss
Research Group FOR916
(Statistical Regularization and Qualitative Constraints) with grant
number 20PA20E-134495/1.}

%
\received{\smonth{1} \syear{2014}}
%
\revised{\smonth{1} \syear{2015}}

%
\begin{abstract}
In this paper, we deal with the data-driven selection of
multidimensional and
possibly anisotropic bandwidths in the general framework of kernel empirical
risk minimization. We propose a universal
selection rule, which leads to optimal adaptive results in a large
variety of
statistical models such as nonparametric robust regression and
statistical learning with
errors in variables. These results are stated in the context of smooth
loss functions, where
the gradient of the risk appears as a good criterion to measure the
performance of our
estimators. The selection rule consists of a comparison of gradient
empirical risks. It can be
viewed as a nontrivial improvement of the \mbox{so-called}
Goldenshluger--Lepski method to nonlinear estimators. Furthermore, one
main advantage of our
selection rule is the nondependency on the
Hessian matrix of the risk, usually involved in standard adaptive procedures.
\end{abstract}

%
\end{frontmatter}

\section{Introduction}\label{sintro}
We consider the minimization problem of an unknown risk function $\R\dvtx \bR^m\to\bR$, where
$m\geq1$ is the dimension of the statistical
model. We assume the existence of a risk minimizer
%
\begin{equation}
\label{oracle} \f^\star\in\arg\min_{\f\in\bR^m}\R(\f),
\end{equation}
where the risk function corresponds to the expectation of an
appropriate loss function w.r.t. an
unknown distribution. In empirical risk minimization, this
quantity is usually estimated by its empirical version
from an i.i.d. sample.
However, in many problems such as local $M$-estimation or
errors-in-variables models, a
nuisance parameter can be involved in the empirical version. This
parameter most often
coincides with some bandwidth related to a kernel that gives rise to ``kernel
empirical risk minimization.'' One
typically deals with this issue in pointwise estimation, as, for
example, in
Polzehl and Spokoiny \cite{PolzehlSpokoiny06} with localized
likelihoods or in Chichignoud and Lederer
\cite{ChichignoudLederer13} with local $M$-estimators. In learning
theory, many authors
have
recently investigated
supervised and unsupervised learning with errors in variables. As a
rule, such matters require one to plug deconvolution
kernels into the empirical risk, as Loustau and Marteau \cite
{pinkfloyds} in noisy
discriminant analysis
or Hall and Lahiri \cite{HallLahiri08} in quantile and moment
estimation; see also Dattner, Rei{ss} and Trabs
\cite{DattnerReissTrabs13}.

In the above papers, the authors studied the theoretical properties of
kernel empirical risk
minimizers and proposed
deterministic
choices of bandwidths to deduce optimal minimax results. As usual,
these optimal
bandwidths are related to the smoothness of the target function or the
underlying density and
are not achievable in practice. Adaptivity is therefore one of the biggest
challenges. In this respect, data-driven bandwidth
selections have been already proposed in
\cite
{ChichignoudLederer13,ChichignoudLoustau13,DattnerReissTrabs13,PolzehlSpokoiny06},
which are all based on Lepski-type
procedures.


Lepski-type
procedures are rather appropriate to construct
data-driven bandwidths involved in kernels; for further details, see,
for example,
\cite{Katkovnik99,Lepski90,LepskiMammenSpokoiny97}. It is
well known that
they suffer from the restriction to isotropic bandwidths with
multidimensional data,
which is the
consideration of nested neighborhoods (hyper-cube). Many improvements
were made by
Kerkyacharian, Lepski and Picard \cite{KerkyacharianLepskiPicard01}
and more recently by Goldenshluger and Lepski
\cite{GoldenshlugerLepski11}
to select anisotropic bandwidths (hyper-rectangle). Nevertheless, their
approach still does not
provide
anisotropic bandwidth selection for nonlinear estimators, which is the
scope of this paper. The only work we
can mention is
\cite{ChichignoudLederer13} in a restrictive case, which is pointwise
estimation in
nonparametric regression.
Therefore, the study of data-driven selection of anisotropic bandwidths
is still an open issue. Moreover, this field is of great
interest in practice, especially in image denoising; see, for example,
\cite{CastroSalmonWillett12,KatkovnikFoiEgiazarianAstola10}.

The main contribution of our paper is to bring new insights to the
problem of bandwidth selection in kernel empirical risk minimization in
a possible anisotropic framework. To this end, we first introduce a new
criterion called \textit{gradient excess risk}, which makes the
anisotropic bandwidth selection possible. We then provide a novel
data-driven selection based on the comparison of ``Gradient empirical
risks.'' That can be viewed as an extension of the so-called
Goldenshluger--Lepski method (GL method; see \cite
{GoldenshlugerLepski11}) and of the empirical risk comparison method
(ERC method; see \cite{ChichignoudLoustau13}). Eventually, we derive
an upper bound for the gradient excess risk (called gradient
inequality) and optimal results in many settings, such as pointwise and
global estimation in
nonparametric regression and clustering with errors in variables.

Note that we consider the risk minimization over the finite dimensional
set $ \bR^m$. In statistical learning or nonparametric
estimation, one usually aims at estimating a functional object
belonging to some Hilbert space.
However,
in many examples, the target function can be approximated by a finite
object, thanks, for instance, to a
suitable decomposition in a
basis of the Hilbert space. This is typically the case in local $M$-estimation,
where the
target function is assumed to be locally polynomial (and even constant
in many cases).
Moreover, in statistical learning, one is often interested in the
estimation of a finite number
of parameters, as in clustering. The
extension to
the infinite-dimensional case is discussed in Section~\ref{sdiscussion}.

The structure of this paper is as follows: the
main ideas behind the gradient excess risk are introduced in the
remainder of this section. An upper bound for the gradient excess risk
of the data-driven procedure is presented in Section~\ref{gradientinequality}.
This procedure is applied to clustering in Section~\ref{sectionkmeans} and to robust
nonparametric regression in Section~\ref{sectionlocalglobal}.
Additionally, a discussion of our assumptions and an outlook are given
in Section~\ref{sdiscussion}, and Section~\ref{ssimu} illustrates
the behavior of the method with numerical results. The proofs are finally
conducted in the \hyperref[appendix]{Appendix}.

\subsection{The gradient excess risk approach}\label{sgradient}
In the literature, such as in statistical learning, the excess risk
$R(\widehat\theta)-R(\theta^\star)$ is the main criterion to
measure the performance of some estimator $\widehat\theta$. Originally,
Vapnik and Chervonenkis \cite{vapnikold}
proposed to control this quantity via the empirical process theory, which
gives rise to slow rates $\cO( n^{-1/2})$ for the excess risk; see
also \cite{vapnik98}.
In the
last decade, many authors have improved such a bound by giving fast
rates $\cO( n^{-1})$ using the so-called
localization technique; see
\cite{svm,kolt,mammen,nedelec,mendelsonkernel,tsybakov2004}
and Boucheron, Bousquet and Lugosi \cite{surveylugosi} for an overview
in classification. This technique consists of studying the increments
of an empirical process in the neighborhood of the target $\f^\star$.
In particular, it requires a variance-risk correspondence, equivalent
to the eminent margin assumption.
As far as we know, this complicated modus operandi is the major
obstacle to the anisotropic bandwidth selection issue. In what follows,
we introduce an alternative criterion to solve this issue, namely the
gradient excess risk ($G$-excess risk, for short, in the sequel). This
quantity is defined as
%
\begin{equation}
\label{dexcessrisk} \bigl\llvert \D\bigl(\fn,\f^\star\bigr)\bigr\rrvert
_2:=\bigl\llvert \D(\fn)-\D\bigl(\f^\star\bigr)\bigr\rrvert
_2\qquad\mbox{where }\D:=\nabla\R,
\end{equation}
whereas \mbox{$ \llvert  \cdot\rrvert  _2 $} denotes the Euclidean norm on $ \bR^m $ and
$\nabla\R\dvtx \bR^m\to\bR^m$
denotes the gradient of the risk $\R$. With a slight abuse of
notation, $G$ denotes the
gradient, whereas $G(\cdot,\f^\star)$ denotes the $\D$-excess risk.
Under regularity assumptions on $R(\cdot)$, the $G$-excess risk is
linked with the excess risk, thanks to the following lemma.

\begin{lemma}
\label{lemmadmargin}
Let $\f^\star$, defined as in (\ref{oracle}), and $U$ be the
Euclidean ball of $\bR^m$ centered at
$\f^\star$, with radius $\delta>0$. Assume $\f\mapsto\R(\f)$ is
$\cC^2(U)$, each second partial derivative of $\R$ is bounded on $U$
by a constant
$\kappa_1$ and the Hessian matrix $H_{\R}(\cdot)$ is positive
definite at $\f^\star$. Then, for $\delta>0$
small enough, we have
\[
\sqrt{\R(\f)-\R\bigl(\f^\star\bigr)}\leq 2\frac{\sqrt{m\kappa_1}}{\lambda_{\min}}\bigl\llvert
\D\bigl(\f,\f^\star\bigr)\bigr\rrvert _2\qquad \forall\f\in U,
\]
where $\lambda_{\min}$ is the smallest eigenvalue of
$H_{\R}(\f^\star)$.
\end{lemma}

The proof is based on the inverse function theorem and a
Taylor expansion of the function $R(\cdot)$.
Let us explain how this lemma, together with standard probabilistic
tools, leads to fast rates for the excess risk.
In this section, $\Rn$ denotes the usual empirical risk with associated
gradient $\Dn:=\nabla\Rn$ and associated ERM $\fn$ for ease of
exposition. Under the assumptions of Lemma \ref{lemmadmargin}, $\D(\f
^\star)=\Dn(\fn)=(0,\dots,0)^\top$, and we
have the following heuristic:
%
\begin{eqnarray}\label{heuristichuber}
\sqrt{\R(\fn)-\R\bigl(\f^\star\bigr)}&\lesssim&\bigl\llvert
\D\bigl(\fn,\f^\star\bigr)\bigr\rrvert _2 = \bigl\llvert \D(
\fn)-\Dn(\fn)\bigr\rrvert _2
\nonumber\\[-8pt]\\[-8pt]\nonumber
&\leq& \sup_{\f\in\bR^m}\bigl
\llvert \D(\f)-\Dn(\f)\bigr\rrvert _2\lesssim n^{-1/2},
\end{eqnarray}
where $\lesssim$ denotes the inequality up to some positive constant.
The last bound only requires a concentration inequality applied to the
empirical process $ \Dn(\cdot)-\D(\cdot) $. Therefore, this
heuristic provides fast rates for the excess risk without any
localization technique. Furthermore, similar bounds can be obtained for
the $\ell_2$-norm $ \llvert  \fn-\f^\star\rrvert  _2 $ using the same path. Indeed,
under the same assumptions, the assertion of Lemma \ref{lemmadmargin}
holds, replacing the square root of the excess risk by $ \llvert  \fn-\f
^\star\rrvert  _2 $ (see the proof of Lemma \ref{lemmadmargin}), and then
optimal rates are deduced.


From the model selection point of view, standard penalization
techniques---based on localization---suffer from the dependency on
parameters involved in the margin assumption. More precisely, in
the strong margin assumption framework, the construction of the penalty
requires the knowledge of
$\lambda_{\min}$, related to the Hessian matrix of the risk. Although
many authors have
recently investigated the adaptivity w.r.t. these parameters, by proposing
``margin-adaptive'' procedures (see \cite{PolzehlSpokoiny06} for the
propagation method, \cite{Lecue07} for aggregation
and \cite{ArlotMassart09} for the
slope heuristic), the theory
is not completed and remains a hard issue; see the related discussion in
Section~\ref{sdiscussion}.
As an alternative, our data-driven procedure
does not suffer from the dependency on $\lambda_{\min}$ since we
focus on a gradient inequality in Section~\ref{gradientinequality}.

\subsection{Kernel empirical risk minimization}\label{sectionKERM}
In this section, the kernel empirical risk minimization is properly
defined and illustrated with two examples: local \mbox{$M$-}estimators and
deconvolution $k$-means.
For some $p\in\bN^\star$, consider a
\mbox{$\bR^p$-}random variable $Z$ distributed according to $P$, absolutely
continuous w.r.t. the
Lebesgue measure. In what follows, we observe a sample
$\cZ_n:=\{Z_1,\ldots,Z_n\}$ of independent
and identically distributed (i.i.d.) random\vspace*{1pt} variables according to $P$.
Moreover, we call a kernel of order $
r\in\bN^\star$ a symmetric function
$K\dvtx \bR^d\to\bR$, $d\geq1$, which satisfies the following properties:
\begin{itemize}
\item[$\bullet$] $\int_{\bR^d} K(x)\,dx=1$,\vspace*{1pt}
\item[$\bullet$] $\int_{\bR^d} K(x)x_j^k\,dx=0$ $\forall k\leq r, \forall j\in\{1,\ldots,
d\}$,\vspace*{1pt}
\item[$\bullet$] $\int_{\bR^d} \llvert  K(x)\rrvert  \llvert  x_j\rrvert  ^{r}\,dx<\infty$,
$\forall j\in\{1,\ldots, d\}$.
\end{itemize}
For any $h\in\cH\subset\bR^d_+$, the dilation
$K_h$ is defined as
\[
K_h(x)=\Pi_h^{-1}K(x_1/h_1,
\ldots, x_d/h_d)\qquad \forall x\in\bR^d,
\]
where $\Pi_h:=\prod_{j=1}^dh_j$.
For a given kernel $K$, we define the kernel
empirical risk
indexed by an anisotropic bandwidth $h\in\cH\subset(0,1]^d$ as
%
\begin{equation}
\label{defemprisk} \Rn_{h}(\f):= \frac{1}{n}\sum
_{i=1}^n\ell_{K_h}(Z_i,\f),
\end{equation}
and an associated kernel empirical risk minimizer (kernel ERM) as
%
\begin{equation}
\label{defkerm} \fn_h\in\argmin_{\f\in\bR^m} \Rn_{h}(
\f).
\end{equation}
The function $\ell_{K_h}\dvtx \bR^p\times\bR^m\to\bR_+$ is a loss function
associated to a kernel $ K_h $ such that $ \f\mapsto\ell_{K_h}(Z,\f
) $ is twice differentiable
$P$-almost surely and such that the limit of its expectation coincides
with the risk, that is,
%
\begin{equation}
\label{eqlimitrisk} \lim_{h\to(0,\dots,0)}\bE\Rn_h(\f)=\R(\f)\qquad \forall\f\in\bR^m,
\end{equation}
where $ \bE$ denotes the expectation w.r.t. the distribution of the
sample $ \cZ_n $.

The agenda is the data-driven selection of the ``best'' estimator in
the
family $\{\fn_h,h\in\cH\}$. This issue arises in many examples, such
as local fitted likelihood (Polzehl and Spokoiny \cite
{PolzehlSpokoiny06}), image denoising
(Astola et~al. \cite{KatkovnikFoiEgiazarianAstola10}) and robust
nonparametric regression; see
Chichignoud and Lederer \cite{ChichignoudLederer13}. In such a
framework, we observe a sample of i.i.d. pairs $Z_i=(W_i,Y_i)_{i=1}^n$,
and the kernel empirical risk has the following general form:
\[
\frac{1}{n}\sum_{i=1}^n
\ell_{K_h}(Z_i,\theta)=\frac{1}{n}\sum
_{i=1}^n\rho(Z_i,\theta){K}_h
(W_i-x_0 ), %
\]
where $\rho(\cdot,\cdot)$ is some likelihood and $x_0\in\bR^d$.
Another example arises when we observe a contaminated sample
$Z_i=X_i+\varepsilon_i$, $i=1, \ldots, n$ in the problem of clustering.
In this case, the kernel empirical risk is defined according to
\[
\frac{1}{n}\sum_{i=1}^n
\ell_{K_h}(Z_i,\bc)=\frac{1}{n}\sum
_{i=1}^n\int_{\bR^d}\min
_{j=1,\dots,k}\llvert x-c_j\rrvert _2^2
\widetilde K_h(Z_i-x)\,dx, %
\]
where $\widetilde K_h(\cdot)$ is a deconvolution kernel and $\bc
=(c_1, \ldots, c_k)\in\bR^{dk}$ is a codebook.

In the next section, we present the bandwidth selection rule in the
general context
of kernel empirical risk minimization. We especially deal with
clustering with errors in variables and robust nonparametric
regression in Sections~\ref{sectionkmeans} and \ref{sectionlocalglobal}, respectively.

\section{Selection rule and gradient inequality}
\label{gradientinequality}
The anisotropic bandwidth selection issue has been recently
investigated in Goldenshluger and Lepski
\cite{GoldenshlugerLepski11} (GL~method)
in density estimation; see also \cite{ComteLacour13} for
deconvolution estimation and
\cite{GoldenshlugerLepski08,GoldenshlugerLepski09} for the white
noise model.
This method, based on the comparison of estimators, requires some ``linearity''
property, which is trivially satisfied by kernel estimators. However,
kernel ERMs are usually
nonlinear (except for the least square estimator), and the GL method
cannot be directly applied to such estimators.

To tackle this issue, we introduce a new selection rule based on the
comparison of
gradient empirical risks instead of estimators (i.e., kernel ERM). To
that end, we first introduce some notations. For any
$h\in\cH$ and any
$\f\in\bR^m$, the gradient empirical risk ($G$-empirical risk) is
defined as
%
\begin{equation}
\label{defdemprisk} \Dn_h(\f):=\frac{1}{n}\sum
_{i=1}^n\nabla\ell_{K_h}(Z_i,
\f)= \Biggl(\frac{1}{n}\sum_{i=1}^n
\frac
{\partial}{\partial\f_j}\ell_{K_h} (Z_i,\f) \Biggr)_{j=1,\ldots, m}.
\end{equation}
Note that we have coarsely $ \Dn_h(\fn_h)=(0,\ldots,0)^\top$ since
$\ell_{K_h}(Z_i,\cdot)$ is
twice differentiable almost surely.
According to (\ref{eqlimitrisk}), we also notice that the limit of
the expectation of the
$\D$-empirical risk coincides with the
gradient of the risk.

Following Goldenshluger and
Lepski \cite{GoldenshlugerLepski11}, we introduce an auxiliary\break \mbox{$\D
$-}empiri\-cal risk in the comparison.
For any couple of bandwidths $(h,\eta)\in\cH^2$ and any $ \f\in\bR
^m $, the auxiliary
$\D$-empirical risk is defined
as
%
\begin{equation}
\label{defconvolutiondemprisk} \Dn_{h,\eta}(\f):=\frac{1}{n}\sum
_{i=1}^n\nabla\ell_{K_h*K_\eta
}(Z_i,
\f),
\end{equation}
where $K_h*K_\eta(\cdot):=\int_{\bR^d}K_h(\cdot-x)K_\eta(x)\,dx$
stands for the convolution
between $K_h$ and $K_\eta$.
The gradient inequality stated in Theorem \ref{thmainresult} is
based on the control of some
random processes as follows.

\begin{defi}[(Majorant)]\label{defmajorant}
For any integer $ {l}>0 $, we call \textit{majorant}
a function $ \majorant_l\dvtx \cH^2\to\bR_+ $ such that
\[
\bP \Bigl(\sup_{\lambda,\eta\in\cH} \bigl\{ \llvert \Dn_ { \lambda,\eta} -\bE
\Dn_{\lambda,\eta}\rrvert _{2,\infty}+\llvert \Dn_{\eta}-\bE
\Dn_{\eta
}\rrvert _{2,\infty}-\majorant_l(\lambda,
\eta) \bigr\}_+>0 \Bigr)\leq n^{-l},
\]
where $\llvert  T\rrvert  _{2,\infty}:=\sup_{\f\in\bR^m}\llvert  T(\f)\rrvert  _2$ for all $
T\dvtx \bR^m\to\bR^m
$ with $\llvert  \cdot\rrvert  _2$ the Euclidean norm on $\bR^m$, and $\bE$ is
understood coordinatewise.
\end{defi}

The main issue for applications is
to compute right order majorants. It follows from classical tools such
as Talagrand's inequalities
(Talagrand \cite{talagrandinventiones},
Boucheron, Lugosi and Massart \cite{boucheronlivre},
Bousquet \cite{bousquet}; see also \cite{GoldenshlugerLepski09b}). In Sections~\ref{sectionkmeans} and~\ref{sectionlocalglobal} such majorant
functions are computed in
clustering and in robust nonparametric regression.

We are now ready to define the selection rule as
%
\begin{equation}
\label{defrule} \widehat h\in\argmin_{h\in\cH} \widehat{\mathrm{BV}}(h),
\end{equation}
where $\widehat{\mathrm{BV}}(h)$ is an estimate of the bias--variance
decomposition at a given
bandwidth $h\in\cH$. It is explicitly defined as
%
\begin{eqnarray}\label{defempbv}
\widehat{\mathrm{BV}}(h):=\sup_{\eta\in\cH} \bigl\{
\llvert \Dn_{h,\eta
}-\Dn_\eta\rrvert _{2,\infty} -
\majorant_l(h,\eta) \bigr\} +\majorant^\infty_l(h)\nonumber
\\
\eqntext{\displaystyle\mbox{with }\majorant_l^\infty (h):=\sup
_{\lambda\in\cH} \majorant_l(\lambda,h).}
\end{eqnarray}
%
The kernel ERM $\fn_{\widehat h}$, defined in (\ref{defkerm}), with
bandwidth $\widehat h$, selected in (\ref{defrule}), satisfies the
following bound.

\begin{theorem}[(Gradient inequality)]\label{thmainresult}
For any $n\in\bN^\star$ and for any $l\in\bN^\star$, we have with
probability $1-n^{-l}$,
\[
\bigl\llvert \D\bigl(\fn_{\widehat
h},\f^\star\bigr)\bigr\rrvert
_2\leq3\inf_{h\in\cH} \bigl\{B(h)+\majorant
_l^\infty(h) \bigr\}, %
\]
where $B\dvtx \cH\to\bR_+$ is a bias function defined as
%
\begin{equation}
\label{defbias} \quad\bias(h):=\max \Bigl(\llvert \bE\Dn_{h}-\D\rrvert
_{2,\infty},\sup_{\eta\in
\cH}\llvert \bE\Dn_{h,\eta}-\bE
\Dn_\eta\rrvert _{
2, \infty} \Bigr)\qquad\forall h\in\cH.
\end{equation}
\end{theorem}

Theorem \ref{thmainresult} is the main result of this
paper. The
$\D$-excess risk of the data-driven estimator $\fn_{\widehat h}$ is
bounded with high probability. The RHS in the gradient inequality can
be viewed as the minimization of a usual
bias--variance trade-off. Indeed, the bias term $ B(h) $ is
deterministic and tends to $0$ as $
h\to(0,\dots,0)$. The majorant $ \majorant_l^\infty(h) $ upper
bounds the stochastic part of
the $ \D
$-empirical risk and can be viewed as a variance term. 

The gradient inequality of Theorem \ref{thmainresult} is sufficient to
establish adaptive fast rates in noisy clustering and
adaptive minimax rates in nonparametric estimation; see Sections~\ref{sectionkmeans}~and~\ref{sectionlocalglobal}. Moreover, the
construction of the selection rule (\ref{defrule}), as
well as the upper bound in Theorem \ref{thmainresult}, does not
suffer from the dependency on
$\lambda_{\min}$ related to the smallest eigenvalue of the Hessian
matrix of the risk; see Lemma \ref{lemmadmargin}. In other words, the
method is robust w.r.t. this
parameter, which is
a major improvement in comparison with other adaptive or model
selection methods
of the literature cited in the \hyperref[sintro]{Introduction}.

\begin{pf*}{Proof of Theorem \protect\ref{thmainresult}}
For some $h\in\cH$, we start with the following decomposition:
%
\begin{eqnarray}
\label{eqdecompositionexcessrisk} \bigl\llvert \D\bigl(\fn_{\widehat h},\f^\star\bigr)
\bigr\rrvert _2&=&\bigl\llvert (\Dn_{\widehat
h}-\D) (\fn_{\widehat h})\bigr\rrvert _2\leq\llvert \Dn_{\widehat
h}-\D \rrvert _{2,\infty}
\nonumber\\[-8pt]\\[-8pt]\nonumber
&\leq&\llvert \Dn_{\widehat h}-\Dn_{\widehat h,h}\rrvert _{2,\infty}+
\llvert \Dn _{\widehat
h,h}-\Dn_{h}\rrvert _{2,\infty}+\llvert
\Dn_{h}-\D\rrvert _{2,\infty}.
\end{eqnarray}
By definition of $\widehat h$ in (\ref{defrule}), the first two terms
in the RHS of (\ref{eqdecompositionexcessrisk}) are bounded as follows:
\begin{eqnarray}\label{eqcontrolvar}
&& \llvert \Dn_{\widehat h}-\Dn_{\widehat h,h}\rrvert
_{2,\infty}+\llvert \Dn_{\widehat h,h}-\Dn_{h}\rrvert
_{2,\infty}\nonumber
\\
&&\qquad = \llvert \Dn_{h,\widehat h}-\Dn_{\widehat
h}\rrvert _{2,\infty}-\majorant_{\ell}(h,\widehat h)+\majorant_{\ell}(h,
\widehat h)
\nonumber
\\
&&\quad\qquad {}+\llvert \Dn_{\widehat h,h}-\Dn_{h}\rrvert _{2,\infty}-
\majorant_{\ell
}(\widehat h,h)+\majorant_{\ell}(\widehat h,h)
\nonumber\\[-8pt]\\[-8pt]\nonumber
&&\qquad \leq \sup_{\eta\in\cH} \bigl\{\llvert \Dn_{h,\eta}-
\Dn_{\eta
}\rrvert _{2,\infty}-\majorant_{\ell}(h, \eta)
\bigr\} +\majorant_ {
\ell}^\infty( h)
\nonumber
\\
&&\quad\qquad{}+ \sup_{\eta\in\cH} \bigl\{\llvert \Dn_{\widehat h,\eta}-
\Dn_{\eta
}\rrvert _{2,\infty}-\majorant_{\ell}(\widehat h,
\eta) \bigr\}+\majorant_{\ell}^\infty(\widehat h)
\nonumber
\\
&&\qquad =\widehat{\mathrm{BV}}(h)+\widehat{\mathrm{BV}}(\widehat h)\leq2\widehat{
\mathrm{BV}}(h).\nonumber
\end{eqnarray}
Besides, the last term in (\ref{eqdecompositionexcessrisk}) is
controlled as follows:
\begin{eqnarray*}
\llvert \Dn_{h}-\D\rrvert _{2,\infty}&\leq&\llvert
\Dn_{h}-\bE\Dn_h\rrvert _{2,\infty}+\llvert \bE
\Dn_{h}-\D\rrvert _{2,\infty}
\nonumber
\\
&\leq&\llvert \Dn_{h}-\bE\Dn_h\rrvert _{2,\infty}-
\majorant_l(\lambda,h) +\majorant_l(\lambda,h)+\llvert
\bE\Dn_{h}-\D\rrvert _{2,\infty}
\nonumber
\\
&\leq&\sup_{\lambda,\eta}
\bigl\{\llvert \Dn_{\lambda,\eta}-\bE\Dn _{\lambda,\eta}\rrvert
_{2,\infty}+\llvert \Dn_{
\eta} -\bE\Dn_\eta\rrvert
_{2,\infty} -\majorant_l(\lambda, \eta) \bigr\}
\nonumber
\\
&&{}+\majorant_l^\infty(h)+\llvert \bE\Dn_{h}-\D
\rrvert _{2,\infty}
\nonumber
\\
&=:&\zeta+\majorant_l^\infty(h)+\llvert \bE
\Dn_{h}-\D\rrvert _{2,\infty}.
\end{eqnarray*}
Using (\ref{eqdecompositionexcessrisk}) and (\ref{eqcontrolvar}), together with the last
inequality, we have for all $
h\in\cH$,
%
\begin{equation}
\label{eqbounddexcessrisk} \bigl\llvert \D\bigl(\fn_{\widehat
h},\f^\star\bigr)
\bigr\rrvert _2\leq2\widehat{\mathrm{BV}}(h)+\zeta+\majorant
_l^\infty(h)+\llvert \bE\Dn_{h}-\D\rrvert
_{2,
\infty}.
\end{equation}
It then remains to control the term $\widehat{\mathrm{BV}}(h)$. We have
\begin{eqnarray*}
&& \widehat{\mathrm{BV}}(h)-\majorant_l^\infty(h)
\\
&&\qquad \leq \sup
_{\lambda,\eta} \bigl\{\llvert \Dn_{\lambda,\eta}-\bE
\Dn_{\lambda,\eta}\rrvert _{2,\infty}+\llvert \Dn_{\eta}-\bE
\Dn_{\eta}\rrvert _{2,
\infty} -\majorant_l(\lambda,
\eta) \bigr\}
\\
&&\quad\qquad{}  +\sup_{\eta}\llvert \bE\Dn_{h,\eta} -\bE
\Dn_\eta\rrvert _{2,\infty}
\\
&&\qquad = \zeta+\sup_{\eta}
\llvert \bE\Dn_{h,\eta} -\bE\Dn_\eta\rrvert _{2,\infty}.
\end{eqnarray*}
The gradient inequality follows directly from (\ref{eqbounddexcessrisk}), Definition \ref{defmajorant} and the definition
of $\zeta$.
\end{pf*}

\section{Application to noisy clustering}\label{sectionkmeans}
Let us
consider an integer $k\geq1$ and a $\bR^d$-random variable $X$ with
law $P$ with density $f$
w.r.t. the Lebesgue measure on $\bR^d$ satisfying $\bE
_P\llvert  X\rrvert  ^2_2<\infty$, where
\mbox{$\llvert  \cdot\rrvert  _2$} stands for the Euclidean norm in $\bR^d$. Moreover, we
restrict the study to the compact set $[0,1]^d$, assuming that $X\in
[0,1]^d$ almost surely. We want to construct $k$
centroids
minimizing some distortion,
%
\begin{equation}
\label{distortion} \RC(\bc):=\bE_{P}w(\bc,X),
\end{equation}
where $\bc=(c_1,\ldots,c_k)\in\bR^{d\times k}$ is a candidate
codebook of $k$ centroids. For ease of
exposition, we study this quantization problem with the Euclidean
distance, by choosing the standard $k$-means loss function, namely,
\[
w(\bc,x)=\min_{j=1,\ldots, k}\llvert x-c_j\rrvert
_2^2,\qquad x\in\bR^d. %
\]
In this section, we are interested in the inverse statistical learning
context (see~\cite{isl}), which corresponds to the minimization of (\ref
{distortion}), thanks to a noisy set
of observations,
\[
Z_i=X_i+\varepsilon_i, \qquad i=1, \ldots, n,
\]
where $(\varepsilon_i)_{i=1}^n$ are i.i.d. with density $g$ w.r.t. the Lebesgue
measure
on $\bR^d$ and mutually independent of the original sample
$(X_i)_{i=1}^n$. This topic was first
considered in \cite{bl12}, where general oracle inequalities are
proposed. Let us fix a
kernel $K_h$ of order $ r\in\bN^\star$ with $h\in\cH$ and consider
$ \widetilde K_h$ a deconvolution kernel
defined such that $\cF[\widetilde
K_h]=\cF[K_h]/\cF[g]$, where $\cF$ stands for the usual Fourier
transform. As introduced in
Section~\ref{sectionKERM}, we have at our disposal the family of
kernel ERM
defined
as
%
\begin{equation}
\label{noisykmeans} \qquad\bcn_h\in\arg\min_{\bc\in\bR^{dk}}
\RCn_h(\bc)\qquad\mbox{where }\RCn_h(\bc):=\frac{1}{n}
\sum_{i=1}^n w(\bc,\cdot)*\widetilde
K_h(Z_i-\cdot),
\end{equation}
where $f*g(\cdot):=\int_{[0,1]^d}f(x)g(\cdot-x)\,dx$ stands for the
convolution product
(restricted to
the compact $[0,1]^d$ for simplicity). From an adaptive point of view,
Chichignoud and Loustau \cite{ChichignoudLoustau13} have recently
investigated the problem of choosing the bandwidth
in (\ref{noisykmeans}). They established fast rates of
convergence---up to a logarithmic term---for a
data-driven selection of $h$, based on a comparison of kernel empirical
risks. However, their approach is restricted to isotropic bandwidth
selection and depends on the parameters involved in the margin
assumption (in
particular on $\lambda_{\min}$
in Lemma \ref{lemmadmargin}).

In the following, adaptive fast rates of convergence for the excess
risk are obtained via the gradient approach. For this purpose, we
assume that the Hessian matrix $H_{\RC}$ is positive definite. This
assumption was considered for the first time in Pollard \cite
{pollard81} and is often referred as Pollard's
regularity assumptions; see also~\cite{levrard}. Under these
assumptions, we can state the same kind of result as Lemma
\ref{lemmadmargin} in the framework of clustering with $k$-means.

\begin{lemma}
\label{dclustering}
Let $\bc^\star$ be a minimizer of (\ref{distortion}), and assume $f$
is continuous and $H_{\RC}(\bc^\star)$ is
positive definite. Let $U$ be the Euclidean ball center at $\bc^\star
$ with radius $\delta>0$. Then, for $\delta$ sufficiently small,
\[
\sqrt{\RC(\bc)-\RC\bigl(\bc^\star\bigr)}\leq C\bigl\llvert \nabla\RC(
\bc)-\nabla\RC\bigl(\bc^\star\bigr)\bigr\rrvert _2\qquad \forall
\bc\in U,
\]
where $C>0$ is a constant which depends on $H_{\RC}(\bc^\star)$, $k$
and $d$.
\end{lemma}
%
We\vspace*{1pt} have at our disposal a family of kernel ERM $\{\bcn_h,h\in\cH\}$
with associated kernel empirical risk
$
\RCn_h(\cdot)
$
defined in
(\ref{noisykmeans}). We
propose to apply the selection rule (\ref{defrule}) to choose the
bandwidth $h\in\cH$. In this
problem as well, we first consider the $\D$-excess risk approach to establish
adaptive fast
rates of convergence for the excess risk. For any $h\in\cH$, the
$ \D$-empirical risk vector of $\bR^{dk}$ is given by
\begin{eqnarray*}
\label{defndemprisk} \DCn_h(\bc)&:=& \Biggl(\frac{1}{n}\sum
_{i=1}^n\frac{\partial
}{\partial c_j^u}\int_{[0,1]^d}w(
\bc,x)\widetilde K_h(Z_i-x)\,dx \Biggr)_{(u,j)\in\{1,
\ldots,d\}\times\{1,\ldots, k\}}
\\
&=& \Biggl(-\frac{1}{n}\sum_{i=1}^n2
\int_{V_j(\bc
)}\bigl(x^u-c_j^u
\bigr)\widetilde K_h(Z_i-x)\,dx \Biggr)_{(u,j)\in\{1, \ldots,d\}\times\{1,\ldots, k\}},
\end{eqnarray*}
where $x^u$ denotes the $u$th coordinate of $x\in\bR^d$ and
$V_j(\bc)$, ${j=1,\ldots, k}$ are open Vorono\"i cells associated to
$\bc$, defined as $V_j(\bc)=\{x\in[0,1]^d\dvtx \forall u\neq j,
\llvert  x-c_j\rrvert  _2<\llvert  x-c_u\rrvert  _2\}$.
Note that $ \DCn_h(\bcn_h)=(0, \ldots, 0)^\top$ a.s. by smoothness.
The construction of the rule follows the general case of Section~\ref{gradientinequality}, which
requires the introduction of an auxiliary $\D$-empirical risk. For any
couple of bandwidths $(h,\eta)\in\cH^2$, the auxiliary
$\D$-empirical risk is defined as
\[
\label{defauxdemprisk} \DCn_{h,\eta}(\bc):= \Biggl(-\frac{1}{n}\sum
_{i=1}^n2\int_{V_j}
\bigl(x^u-c_j^u\bigr)\widetilde
K_{h,\eta}(Z_i-x)\,dx \Biggr)_{(u,j)\in\{1, \ldots,d\}\times\{
1,\ldots, k\}}\in
\bR^{dk},
\]
where $\widetilde K_{h,\eta}=\widetilde{K_h*K_\eta}$ is the
auxiliary deconvolution kernel as in
Comte and Lacour \cite{ComteLacour13}.

The statement of the oracle inequality is based on the computation of a majorant
function. For this purpose, we need the following additional
assumptions on the kernel ${K}\in\bL_2(\bR^d)$.

(\textbf{K1}) There exists $S=(S_1,\dots,S_d)\in\bR^d_+$
such that the kernel
$K$
satisfies
\[
\operatorname{supp}\mathcal{F}[K]\subset[-S,S]\quad\mbox{and}\quad\sup_{t\in\bR
^d}
\bigl\llvert \mathcal{F}[K](t)\bigr\rrvert < \infty, %
\]
where $\operatorname{supp} g=\{x\dvtx g(x)\neq0\}$ and $[-S,S]=\bigotimes_{v=1}^d [-S_v,S_v]$.

This assumption is standard in deconvolution estimation and
is satisfied by
many standard kernels, such as the \textit{sinc} kernel.

We also consider a kernel $K$ of order $r\in\bN^\star$, according to
the definition of Section~\ref{sectionKERM}. 
Kernels of order
$r$ satisfying (\textbf{K1}) could be constructed by using the
so-called Meyer wavelet; see \cite{Mallat09}. Additionally, we need an
assumption on the noise distribution $g$:

\begin{longlist}
\item[\textbf{Noise assumption} \textbf{NA}$(\rho,\beta)$.] There
exist a vector $ \beta=(\beta_1,\dots,\beta_d)\in (0,\infty)^d $
and a positive constant $ \rho$ such that for all $ t\in\bR^d $,
\[
\bigl\llvert \cF[g](t)\bigr\rrvert \geq\rho\prod_{j=1}^d
\biggl(\frac
{t_j^2+1}{2} \biggr)^{-\beta_j/2}. %
\]
\end{longlist}

\textbf{NA($\rho,\beta$)} deals with a polynomial behavior of the
Fourier transform of the noise density $g$. An exponential decreasing
of the characteristic function of $g$ is not considered in this paper
for simplicity; see
\cite{ComteLacour13} in multivariate deconvolution for such a study.

We are now ready to compute some majorant functions in our context.
For some $ s^+>0 $, let $ \cH:=[h_-,h^+]^d $ be the bandwidth set such
that $0<h_-<h^+<1$,
%
\begin{equation}
\label{defhminhmax} h_{-}:= \biggl(\frac{\log^6(n)}{n} \biggr)^{1/(2\vee2\sum
_{j=1}^d\beta_j)}
\quad\mbox{and}\quad h^{+}:= \bigl(1/\log(n) \bigr)^{1/(2s^{+})}.
\end{equation}

\begin{lemma}
\label{lemmamajkmeans}
Assume $(\mathbf{K1})$ and $\mathbf{NA}(\rho,\beta)$ hold for some $\rho
>0$ and some
$\beta\in\bR^d_+$.
Let $a\in(0,1)$, and consider $
\cH_a:=\{(h_{-},\dots,h_{-})\}\cup \{h\in\cH\dvtx  \forall
j=1,\ldots,d\ \exists
m_j\in\bN\dvtx  h_j=h^{+}a^{m_j} \}$ an exponential net of $
\cH=[h_-,h^+]^d $,
such that $ \llvert  \cH_a\rrvert  \leq n $.
For any
integer $ {l}>0 $, let us introduce the function $ \mathcal
{M}^{\mathrm{k}}
_l\dvtx \cH^2\to\bR_+ $
defined
as
\[
\label{defmajorantkmeans} \mathcal{M}^{\mathrm{k}}_l(h,\eta):=
b'_1\sqrt{kd} \biggl(\frac{\prod_{i=1}^d\eta_i^{-\beta_i}}{\sqrt
{n}}+
\frac{\prod_{i=1}
^d(h_i\vee\eta_i)^{-\beta_i}}{\sqrt{n}} \biggr),
\]
where $b'_1:=b'_1(l)>0$ is linear in $l$ and independent of $n$; see
the \hyperref[appendix]{Appendix} for details.

Then, for $n$ sufficiently large, the function $\mathcal{M}^{\mathrm{k}}
_l(\cdot,\cdot)$ is a majorant, that is,
\begin{eqnarray*}
&& \bP \Bigl(\sup_{h,\eta\in\cH_a} \bigl\{ \llvert \DCn_ { h,\eta} -
\bE\DCn_{h,\eta}\rrvert _{2,\infty}+\llvert \DCn_{\eta}-\bE
\DCn_{\eta }\rrvert _{2,\infty}-\mathcal{M}^{\mathrm{k}}_l(h,
\eta) \bigr\}_+>0 \Bigr)
\\
&&\qquad \leq n^{-l},
\end{eqnarray*}
where\vspace*{1pt} $ \bE$ denotes the expectation w.r.t. to the sample and
$\llvert  T\rrvert  _{2,\infty}=\break \sup_{\bc\in[0,1]^{dk}}\llvert  T(\bc)\rrvert  _2$ for all $ T\dvtx \bR
^{dk}\to\bR^{dk}
$ with $\llvert  \cdot\rrvert  _2$ the Euclidean norm on $\bR^{dk}$.
\end{lemma}

The proof is based on a Talagrand inequality; see
the \hyperref[appendix]{Appendix}. This lemma is
the
cornerstone and gives the order of the variance term in such a
problem.


We are now ready to define the selection rule in this setting as
%
\begin{equation}
\label{defrulekmeans} \qquad\widehat h\in\argmin_{h\in\cH_a} \Bigl\{\sup
_{\eta\in\cH_a} \bigl\{\llvert \DCn_{h,\eta}-\DCn_\eta
\rrvert _{2,\infty} -\mathcal{M}^{\mathrm{k}}_l(h,\eta) \bigr
\} +\cM^{\mathrm{k},\infty}_l(h) \Bigr\},
\end{equation}
where
$\cM^{\mathrm{k},\infty}_l(h):=\sup_{\lambda\in\cH_a}
\mathcal{M}^{\mathrm{k}}_l(\lambda,h)
$ and $\cH_a$ is defined in Lemma \ref{lemmamajkmeans}.
Eventually, we need an additional
assumption on the regularity of the density $f$ to control the bias
term in Theorem \ref{thmkmeans}. The regularity is expressed in terms
of anisotropic Nikol'skii class.

\begin{defi}[(Anisotropic Nikol'skii space)]\label{defnikolskiiAnisot}
Let $ s=( s_1, s_2,\ldots, s_d)\in\bR^d_+ $, $ q\in[1,\infty[ $
and $ L>0
$ be fixed.
We say that $ f\dvtx [0,1]^d\rightarrow[-L,L] $ belongs to the anisotropic
Nikol'skii class
$\cN_{q,d}(s,L)$ if for all $ j=1,\ldots,d $, $z\in\bR$ and for all $
x\in(0,1]^d $,
\begin{eqnarray*}
&& \biggl(\int\biggl\llvert \frac{\partial^{\lfloor s_j\rfloor}}{\partial
x_j^{\lfloor s_j\rfloor}}f(x_1,
\ldots,x_j+z,\ldots,x_d)-\frac
{\partial^{\lfloor s_j\rfloor}}{\partial x_j^{\lfloor s_j\rfloor
}}f(x_1,
\ldots,x_j,\ldots,x_d)\biggr\rrvert ^q\,dx
\biggr)^{
1/q}
\\
&&\qquad  \leq L\llvert z\rrvert ^{ s_j-\lfloor s_j\rfloor}, %
\end{eqnarray*}
and $\llVert   \frac{\partial^{ l}}{\partial x_j^{ l}}f\rrVert  _q\leq L$, for any
$l=0,\ldots, \lfloor s_j\rfloor$, where $\lfloor s_j\rfloor$ is the
largest integer strictly less than $s_j$.
\end{defi}

Nikol'skii classes were introduced in approximation theory
by Nikol'skii; see \cite{Nikolskii75}, for example. We also refer to
\cite{GoldenshlugerLepski11,KerkyacharianLepskiPicard01} where the
problem of adaptive
estimation has been treated for
the Gaussian white noise model and for density estimation, respectively.


In the sequel, we assume that the multivariate density $f$ belongs
to the
anisotropic Nikol'skii class $\cN_{2,d}(s,L)$, for some $s\in\bR
_+^d$ and some $L>0$. In other words, the density has possible
different regularities in all directions. The statement of a
nonadaptive upper bound for the excess risk in the anisotropic case
has been already investigated in
\cite{ChichignoudLoustau13}. In the following theorem, we propose the
adaptive version of the previous
cited result, where the bandwidth $\widehat h$ is chosen via the
selection rule (\ref{defrulekmeans}).

\begin{theorem}
\label{thmkmeans}
Assume $ (\mathbf{K1}) $ and $\mathbf{NA}(\rho,\beta)$ hold for some $\rho
>0$ and some
$\beta\in\bR^d_+$. Assume the Hessian matrix of $\RC$ is positive
definite for any $\bc^\star\in\cM$. Then, for any $ s\in(0,s^+]^d
$, any $L>0$, we have
\begin{eqnarray*}
&&\limsup_{n\to\infty} n^{1/(1+\sum_{j=1}^d\beta_j/s_j)}\sup_{f\in
\cN_{2,d}(s,L)}
\bigl[\bE\RC(\bcn_{\widehat h})-\RC\bigl(\bc^\star \bigr) \bigr]<
\infty,
\end{eqnarray*}
where $\widehat h$ is driven in (\ref{defrulekmeans}).
\end{theorem}

This theorem is a direct application of Theorem \ref{thmainresult},
Lemma \ref{dclustering} and the majorant construction.
It gives adaptive fast rates of convergence for the excess risk of~$\bcn_{\widehat
h}$ and significantly improves
the result stated in \cite{ChichignoudLoustau13} for two reasons:
first, the selection rule allows the extension to the
anisotropic case; besides, there is no logarithmic term in the adaptive
rate. In our opinion, the localization technique used in \cite
{ChichignoudLoustau13} seems to be the major obstacle to avoid the
extra $\log n$ term.

%

\section{Application to robust nonparametric regression}\label{sectionlocalglobal}
In this section, we apply the gradient inequality to the framework of
local $M$-estimation in nonparametric robust regression. It will give
adaptive minimax results for nonlinear estimators for both pointwise
and global estimation.

Let us specify the model beforehand. For some $n\in\bN^\star$, we
observe a training set
$\cZ_{n}:=\{(W_i,Y_i), i=1,\ldots, n\}$ of i.i.d. pairs, distributed
according to the probability
measure $ P $ on $[0,1]^d\times\bR$ satisfying the set of equations
%
\begin{equation}
\label{model} Y_i=f^\star(W_i)+
\xi_i,\qquad i=1,\ldots, n.
\end{equation}
We aim at estimating the target function
$f^\star\dvtx [0,1]^d\rightarrow[-B,B]$, $B>0$. The
noise variables
$(\xi_i)_{i\in{1,\ldots,n}}$ are assumed to be i.i.d. with symmetric
density $g_\xi$
w.r.t. the
Lebesgue measure. We also assume $g_\xi$ is continuous at $ 0 $ and $
g_\xi(0)>0 $.
For simplicity, the design points $(W_i)_{i=1}^n$ are assumed to be
i.i.d. according to the uniform law on $[0,1]^d$ (extension to a more
general design is
straightforward), and we assume that
$(W_i)_{i=1}^n$ and $(\xi_i)_{i=1}^n$ are mutually
independent for ease of exposition. Eventually, we restrict the
estimation of $f^\star$
to the closed set $ \cT\subset[0,1]^d $ to avoid discussion on
boundary effects. We will
consider a point $ x_0 \in\cT$ for pointwise estimation and the
$\bL_q(\cT)$-risk for global estimation, for $q\in[1,+\infty)$.

Next, we introduce an estimate of $f^\star(x_0)$ at any $x_0\in\cT$
with the local constant
approach (LCA).
The key idea of LCA, as described, for example, in \cite{Tsybakov08}, Chapter~1, is to
approximate the target function by a constant in a neighborhood of size $
h\in(0,1)^d $ of a given point $x_0$, which corresponds to a model of dimension
$ m=1 $. To deal with heavy-tailed
noises, we especially employ the Huber loss (see \cite{Huber64})
defined as follows. For any scale $\gamma>0$ and $z\in\bR$,
\[
\label{defhubercontrast} \rho_\gamma(z):=\cases{ z^2/2, &\quad if $
\llvert z\rrvert \leq\gamma$,
\cr
\gamma\bigl(\llvert z\rrvert -\gamma/2\bigr), &
\quad otherwise.}
\]
%
The parameter $ \gamma$ selects the level of robustness of the Huber
loss between
the square loss (large value of $ \gamma$) and the absolute loss
(small value of $\gamma$).
Let $ \cH:=[h_-,h^+]^d $ be the
bandwidth set such that $0<h_-<h^+<1$,
\[
\label{defbandwidthnet} h_-:=\frac{\log^{6/d}(n)}{n^{1/d}}\quad\mbox{and}\quad h^+:=
\frac
{1}{\log^2(n)}.
\]
For any $x_0\in\cT$, the
local
{estimator} $
\estim_{h}(x_0) $ of $ f^\star(x_0) $
is defined as
%
\begin{equation}
\label{deflocalestimate} \estim_{h}(x_0):=\argmin_{t\in[- B, B]}
\widehat{R}^{\mathrm{loc}}_{h}(t),\qquad  h\in\cH,
\end{equation}
where $ \widehat{R}^{\mathrm{loc}}_{h}(\cdot):= \frac{1}{n}\sum_{i=1}^n\rho_\gamma(Y_i-\cdot
) K_h(W_i-x_0)$ is the local
empirical risk, and $ K_h $ is a $1$-Lipschitz kernel of order $1$.
We notice that the local empirical risk estimates the local risk $ \R
^{\mathrm{loc}}
(\cdot):=
\bE_{Y\mid W=x_0}\rho_\gamma(Y-\cdot) $ whose $ f^\star(x_0) $ is its
unique minimizer.

In
nonparametric estimation, one is usually interested in pointwise or
global risk instead of
excess risk.
Since Theorem \ref{thmainresult} controls the $ \D$-excess risk of
the adaptive estimator,
we present the following lemma that links the pointwise risk with the
$\D$-excess risk.

\begin{lemma}\label{lemlocalmargincondition}
Assume that
$
\sup_{h\in\cH}\llvert  \estim_{h}(x_0)-f^\star(x_0)\rrvert  \leq\bE\rho_\gamma
''(\xi_1)/4
$ holds.
Then, for all $ h\in\cH$,
\[
\bigl\llvert \estim_{h}(x_0)-f^\star(x_0)
\bigr\rrvert \leq\frac{2}{\bE\rho_\gamma''(\xi
_1)}\bigl\llvert { G^{\mathrm{loc}} \bigl(
\estim_{h}(x_0) \bigr)-G^{\mathrm{loc}}
\bigl(f^\star(x_0) \bigr)}\bigr\rrvert, %
\]
where $ G^{\mathrm{loc}}$ (and, resp., $ \rho_\gamma'' $) denotes the
derivative of
$\R^{\mathrm{loc}}$ (resp., the second derivative of $
\rho_\gamma$).
\end{lemma}

The proof is given in the \hyperref[appendix]{Appendix}. We can also
deduce the same inequality with the $ \bL_q(\cT) $-norm. The\vspace*{2pt}
assumption $
\sup_{h\in\cH}\llvert  \estim_{h}(x_0)-f^\star(x_0)\rrvert  \leq\bE\rho_\gamma
''(\xi_1)/4$ is necessary to use
the theory of differential calculus and can be satisfied by using the
consistency of $
\estim_{h}
$. In this respect, the definitions of $ h_- $ and $ h^+ $ above\vspace*{2pt} imply the
consistency
of all estimators $\estim_{h}, h\in\cH$; for further details, see
below as well as~\cite{ChichignoudLederer13}, Theorem 1.

\subsection{The selection rule in pointwise estimation}\label{sectionlocal}
We now present the application of the selection rule for pointwise
estimation. To compute the procedure, we
define the $ \D$-empirical risk as
%
\begin{equation}
\label{deflocalempiricalderivative} \widehat{G}^{\mathrm{loc}}_{h}(t):=
\frac{\partial\widehat
{R}^{\mathrm{loc}}_{h}}{\partial t}(t) =-\frac{1}{n}\sum_{i=1}^n
\rho_\gamma' (Y_i-t ) K_h(W_i-x_0).
\end{equation}
For two bandwidths $ h,\lambda$, we introduce the auxiliary $ \D
$-empirical risk as
\[
\widehat{G}^{\mathrm{loc}}_{h,\eta}(t):=-\frac{1}{n}\sum
_{i=1}^n\rho_\gamma'
(Y_i-t ) K_{h,\eta
}(W_i-x_0),
\]
where $ K_{h,\eta}:=K_{h}*K_{\eta}$, as before.

To apply the results of Section~\ref{gradientinequality}, we need to
compute optimal majorants of the associated empirical processes.
The construction of such bounds for the pointwise case has already
received attention in the literature; see \cite{ChichignoudLederer13}, Proposition~2. For any integer $l\in\bN^\star$, let us
introduce the function
$\mathcal{M}^{\mathrm{loc}}_l\dvtx \cH^2\to\bR_+$ defined as
\begin{eqnarray*}
&&\mathcal{M}^{\mathrm{loc}}_l(h,\eta):=C_{0}\llVert K
\rrVert _2\sqrt{\bE\bigl[\rho _\gamma'(
\xi _1)\bigr]^2} \biggl(\sqrt{
\frac{l
\log(n)}{n\prod_{j=1}^dh_j\vee\eta_j}}+\sqrt{\frac{l
\log(n)}{n\prod_{j=1}^d\eta_j}} \biggr),
\end{eqnarray*}
where $C_0>0$ is an absolute constant which does not depend on the
model. Then if we set $
\cH_a:=\{(h_{-},\dots,h_{-})\}\cup \{h\in\cH\dvtx  \forall
j=1,\ldots,d\ \exists
m_j\in\bN\dvtx  h_j=h^{+}a^{m_j} \}$, $a\in(0,1)$, an\vspace*{1pt} exponential
net of $
\cH=[h_-,h^+]^d $,
such that $ \llvert  \cH_a\rrvert  \leq n $, for any $ l>0 $, the
function $
\mathcal{M}^{\mathrm{loc}}_l(\cdot,\cdot)$ is a majorant according
to Definition
\ref{defmajorant}.

Eventually, we introduce the data-driven bandwidth following the schema
of the selection rule
in Section~\ref{gradientinequality},
%
\begin{equation}
\label{defruleloc} \qquad\widehat h^{\mathrm{loc}}\in\argmin_{h\in\cH_a} \Bigl\{\sup
_{\eta\in\cH
_a} \bigl\{\bigl\llvert \widehat{G}^{\mathrm{loc}}_{h,
\eta}
-\widehat{G}^{\mathrm{loc}}_\eta\bigr\rrvert _{\infty} -
\mathcal{M}^{\mathrm{loc}}_l(h,\eta) \bigr\} +\cM^{\mathrm{loc},\infty}_l(h)
\Bigr\},
\end{equation}
where $\cM^{\mathrm{loc},\infty}_l(h):=\sup_{h'\in\cH
_a}\mathcal{M}^{\mathrm{loc}}_l(h',h)$. To derive minimax adaptive
rates for local
estimation, we start with the definition of
the anisotropic H\"older class.

\begin{defi}[(Anisotropic H\"{o}lder class)]\label{defholderAnisot}
Let $ s=(s_1,s_2,\ldots,s_d)\in\bR_+^d $ and $ L>0 $ be fixed.
We say that $ f\dvtx [0,1]^d\rightarrow[-L,L] $ belongs to the anisotropic
H\"{o}lder class
$\Sigma(s,L)$ of functions if for all $ j=1,\ldots,d $ and for all $
x\in(0,1]^d $,
\begin{eqnarray*}
&& \biggl\llvert \frac{\partial^{\lfloor s_j\rfloor}}{\partial x_j^{\lfloor
s_j\rfloor}}f(x_1, \ldots, x_j+z,
\ldots, x_d)-\frac{\partial
^{\lfloor s_j\rfloor}}{\partial x_j^{\lfloor s_j\rfloor}}f(x_1, \ldots,
x_j,\ldots, x_d)\biggr\rrvert
\\
&&\qquad \leq L \llvert z\rrvert
^{s_j-\lfloor s_j\rfloor
}\qquad \forall z\in\bR, %
\end{eqnarray*}
and
\[
\sup_{x\in[0,1]^d}\biggl\llvert \frac{\partial^{l}}{\partial
x_j^{l}}f(x)\biggr\rrvert
\leq L\qquad \forall l=0, \ldots, \lfloor s_j\rfloor, %
\]
where $\lfloor s_j\rfloor$
is the largest integer strictly less than $s_j$.
\end{defi}
%

\begin{theorem}\label{thholderadapation} For any $ s\in(0,1]^d $,
any $ L>0 $ and any
$
q\geq1 $, it
holds for all $ x_0\in\cT$,
\[
\limsup_{n\to\infty} \bigl({n}/\log(n) \bigr)^{q\bar s/(2\bar
s+1)}\sup
_{f\in\Sigma({
s},L)}\bE\bigl\llvert \estim_{\widehat
h^\mathrm{loc}}(x_0)-f^\star(x_0)
\bigr\rrvert ^q<\infty,
\]
where $ \bar s:= (\sum_{j=1}^d s_j^{-1} )^{-1} $ denotes
the harmonic
average.
\end{theorem}

The proposed estimator $ \estim_{\widehat h} $ is then
adaptive minimax over
anisotropic H\"{o}lder classes in
pointwise estimation. The minimax optimality of this rate [with the $
\log(n) $ factor] has
been stated by
\cite{Klutchnikoff05} in the white noise model for pointwise
estimation; see also
\cite{GoldenshlugerLepski08}.
For simplicity, we did not study the case of locally polynomial
functions [i.e., $s\in(0,\infty)^d $].

Chichignoud and Lederer \cite{ChichignoudLederer13}, Theorem~2, have shown that the variance of local
$M$-estimators is of order $
\bE[\rho_\gamma'(\xi_1)]^2/n(\bE\rho_\gamma''(\xi_1))^2
$,
and therefore their Lepski-type procedure depends on this quantity.
Thanks to the gradient approach, we obtain the same
result
without the dependency on the parameter $ \bE\rho_\gamma''(\xi_1)
$, which corresponds to $
\lambda_{\min} $ in the general setting. The selection rule
is therefore robust w.r.t. to the fluctuations of this parameter, in
particular when $ \gamma$
is small (median estimator).

\subsection{The selection rule in global estimation}\label{sectionglobal}
The aim of this section is to derive adaptive minimax results for
$\estim_h$ for the $\bL_q$-risk. To this end, we need to modify the
selection rule
(\ref{defruleloc}) including a global ($\bL_q$-norm) comparison
of $ \D$-empirical risks. For this purpose, for all $ t\in\bR$, we
denote the $\D$-empirical risks at a given point $ x_0\in\cT
$ as
\[
\widehat{G}^{\mathrm{loc}}_{h}(t,x_0)=-\frac{1}{n}
\sum_{i=1}^n\rho _\gamma'
(Y_i-t ) K_h(W_i-x_0)
\]
and
\[
\widehat{G}^{\mathrm{loc}}_{h,\eta}(t, x_0)=-
\frac{1}{n}\sum_{i=1}^n
\rho_\gamma' (Y_i-t ) K_{h,\eta}(W_i-x_0),
\]
where the dependence in $x_0$ is explicitly written. Then we define,
for $
q\in[1,\infty[ $ and for any function $ \omega\dvtx \bR\times\cT\to
\bR$, the $ \bL_q $-norm and $
\bL_{q,\infty} $-semi-norm
\[
\bigl\llVert \omega(t,\cdot)\bigr\rrVert _q:= \biggl(\int
_\cT\bigl\llvert \omega(t,x)\bigr\rrvert ^q\,dx
\biggr)^{1/q}\quad\mbox{and}\quad \llVert \omega\rrVert
_{q,\infty}:=\sup_{t\in[-B,B]}\bigl\llVert \omega(t,\cdot)\bigr
\rrVert _q. %
\]

The construction of majorants is based on uniform bounds for $\bL
_q$-norms of empirical
processes. Recently, Goldenshluger and
Lepski investigated this topic
\cite{GoldenshlugerLepski09b}, Theorem 2. For any integer $l\in\bN
^\star$, let us introduce the function
$\Gamma_{l,q}\dvtx \cH\to\bR_+$ defined as
\begin{eqnarray*}
\Gamma_{l,q}(h)&:=& C_{q}\bigl\llVert \rho_\gamma'
\bigr\rrVert _\infty\sqrt{1+l}
\\
&&{} \times
\cases{\displaystyle 4\llVert K\rrVert
_q{\Biggl(n\prod_{j=1}^dh_j
\Biggr)^{-(q-1)/q}}, &\quad if $q\in[1,2[$,
\cr
\displaystyle\frac{30q}{\log(q)}\bigl(
\llVert K\rrVert _2\vee\llVert K\rrVert _q\bigr){\Biggl(n
\prod_{j=1}^dh_j
\Biggr)^{-1/2}}, &\quad if $q\in[2,\infty[$,}
\end{eqnarray*}
where $C_q>0$ is an absolute constant which does not depend on
$n$. Then, for any $ l>0 $, the
function $\mathcal{M}^{\mathrm{glo}}_{l,q}(\lambda,\eta):=\Gamma_{l,q}
(\lambda\vee\eta)+\Gamma_{l,q}(\eta) $
is a majorant according to Definition \ref{defmajorant}.

We finally select the bandwidth according to
\[
\widehat h^{\mathrm{glo}}_q\in\argmin_{h\in\cH} \Bigl\{\sup
_{\eta\in\cH
} \bigl\{\bigl\llVert \widehat{G}^{\mathrm{loc}}_{h,\eta}
-\widehat{G}^{\mathrm{loc}}_\eta\bigr\rrVert _ { q,
\infty} -
\mathcal{M}^{\mathrm{glo}}_{l,q}(h,\eta) \bigr\} +2\Gamma_{l,q}(h)
\Bigr\}.
\]

The above choice of the bandwidth leads to the estimator $ \estim
_{\widehat
h_q^{\mathrm{glo}}} $ with the following adaptive minimax properties
for the $ \bL_q $-risk
over anisotropic
Nikol'skii
classes; see Definition \ref{defnikolskiiAnisot}.

\begin{theorem}\label{thnikolskiiadapation}
For any $
s\in(0,1]^d $, any $ L>0 $ and
any $
q\geq1 $, it
holds that
\[
\limsup_{n\to\infty}\psi_{n,q}^{-1}(s)\sup
_{f\in\cN_{q,d}({
s},L)}\bE\llVert \estim_{\widehat
h_q^{\mathrm{glo}}}-f\rrVert
_q^q<\infty,
\]
where $ \bar s:= (\sum_{j=1}^d s_j^{-1} )^{-1} $ denotes
the harmonic
average and
\[
\psi_{n,q}(s)=\cases{ (1/n )^{q (q-1)\bar s/(q\bar s+q-1)}, &\quad if $q\in [1,2[$,
\cr
(1/{n} )^{q\bar s/(2\bar s+1)}, &\quad if $q\geq2$.}
\]
\end{theorem}

We refer to \cite
{HasminskiiIbragimow90,HasminskiiIbragimov81} for the
minimax optimality
of these rates over Nikol'skii classes. The proposed estimate $\estim
_{\widehat
h^{\mathrm{glo}}_q}$
is then adaptive
minimax. To the best of our knowledge, the minimax adaptivity over
anisotropic Nikol'skii
classes
has never been studied
in regression with possible heavy-tailed noises. We finally refer to
the remarks after
Theorem \ref{thholderadapation}.

\section{Discussion}
\label{sdiscussion}

Our paper solves the general bandwidth selection issue in kernel ERM by
using a novel
selection rule, based on the minimization of an estimate of the
bias--variance decomposition of
the gradient excess risk. This new criterion simultaneously upper
bounds the
estimation error ($ \ell_2 $-norm) and the prediction error (excess
risk) with optimal rates.

One of the key messages we would like to highlight is the following: if
we consider smooth loss functions
and
a family of consistent ERM, fast rates of
convergence are automatically reached,
provided that the Hessian matrix of the risk function is positive
definite. This statement is
based on the key Lemma \ref{lemmadmargin} in Section~\ref{sgradient}, where the square root of
the excess risk is controlled by the $\D$-excess risk.

From an adaptive point of view, one can take another look at Lemma \ref
{lemmadmargin}. On the
RHS of Lemma \ref{lemmadmargin}, the $\D$-excess risk is multiplied
by the constant
$\lambda_{\min}^{-1}$, that is, the smallest eigenvalue of the
Hessian matrix at $\f^\star$. This
parameter is also involved in the margin assumption. As a
result, our selection rule does not depend on this parameter since the
margin assumption is
not required to obtain slow rates for the $\D$-excess risk.
This fact partially solves an issue highlighted by Massart \cite{Massart07}, Section~8.5.2, in the
model selection framework:

\begin{quotation}
\textit{It is indeed a really hard work in this context to design
margin adaptive
penalties.
Of course recent works on the topic, involving
local Rademacher penalties, for instance, provide at least some theoretical
solution to the problem but still if one carefully looks at the
penalties which
are proposed in these works, they systematically involve constants
which are
typically unknown. In some cases, these constants are absolute constants
which should nevertheless considered as unknown just because the numerical
values coming from the theory are obviously over pessimistic. In some other
cases, it is even worse since they also depend on nuisance parameters related
to the unknown distribution.}
\end{quotation}
%
In Section~\ref{ssimu} below, we also illustrate the robustness of
the method with numerical results. An interesting and challenging open
problem would be to employ the gradient approach in the
model
selection framework in order to propose a more robust penalization
technique (i.e., which
does not depend on the parameter $ \lambda_{\min} $).

The gradient approach
requires two main ingredients: the first
one concerns the smoothness of the loss function in terms of
differentiability; the second one affects the dimension of the
statistical model that we have
at hand, which has to be parametric, that is, of finite dimension $m\in
\bN^\star$. From our point of view, the smoothness of the loss
function is not a restriction, since modern
algorithms are usually based---in order to reduce computational
complexity---on some kind of
gradient descent methods in practice. On the
other hand, the
second ingredient might be more restrictive from the model selection
point of view.
An interesting open problem would be to employ the same path when the
dimension $m\geq1$ is possibly larger than $n$, that is, in a
high-dimensional setting.

\section{Numerical results}
\label{ssimu}
For completeness, we illustrate the performance of our selection rule
in the context of clustering with errors in variables, and compare it
to the most recent bandwidth selection procedure in that framework: ERC
method, recently evolved in \cite{ChichignoudLoustau13}. This method
has both theoretical and computational advantages (see also \cite
{KatkovnikSpokoiny08}); however, it only provides isotropic bandwidth
selection. For this reason, our anisotropic selection rule may
outperform ERC method.


The computation of the selection rule (\ref{defrulekmeans}) requires
many optimization steps. We first compute a family of codebooks $\{
\widehat\bc_h,h\in\mathcal{H}\}$ according to (\ref{noisykmeans}),
by using a noisy version of the vanilla $k$-means algorithm. This
technique gives an approximation of the optimal solution (\ref
{noisykmeans}) thanks to an iterative procedure based on Newton
optimization. More theoretical foundations are detailed in \cite
{bl12}. Second, we use parallel execution in order to compute the
comparison of gradient empirical risks.

\subsection*{Experiments}
We generate an i.i.d. noisy sample $\cD_n=\{Z_1,\ldots, Z_n\}$ such
that for any $i=1,\ldots, n$,
%
\begin{eqnarray}
\label{modexp} Z_i=\cases{ X_i^{(1)}+
\varepsilon_i(u), &\quad if $Y_i=1$,
\cr
X_i^{(2)}+\varepsilon_i(u), &\quad if
$Y_i=2$,}
\end{eqnarray}
where\vspace*{1pt} $(X_i^{(1)})_{i=1}^n$ [resp., $(X_i^{(2)})_{i=1}^n$] are i.i.d.
Gaussian with density $f_{\mathcal{N} (0_2,I_2 )}$ (resp.,
$f_{\mathcal{N} ((5,0)^T,I_2 )}$) and $(Y_i)_{i=1}^n$ are
i.i.d. such that $\bP(Y_i=1)=\bP(Y_i=2)=1/2$.
Here, $(\varepsilon_i(u))_{i=1}^n$ are i.i.d. with Gaussian noise with
zero mean $(0,0)^T$ and covariance matrix $\Sigma(u)={1\ \ 0\choose
0\ \ u}$ for $u\in\{1,\ldots, 10\}$. In\vspace*{1pt} this setting, we compare both
adaptive procedures [our selection rule (\ref{noisykmeans}) and ERC
method] to the standard $k$-means with Lloyd's algorithm by computing
the empirical clustering error according to
%
\begin{eqnarray}
\label{cerror} \mathcal{I}_n(\widehat c_1,\widehat
c_2):=\min_{\widehat\bc
=(\widehat c_1,\widehat
c_2),(\widehat c_2,\widehat c_1)}\frac{1}{n}\sum
_{i=1}^n\ind\bigl(Y_i\neq
f_{\widehat
\bc}(X_i)\bigr),
\end{eqnarray}
where $f_{\widehat\bc}(x)\in\arg\min_{j=1, 2}\llvert  x-\widehat\bc
_j\rrvert  ^2_2$ and
$Y_i\in\{1,2\}$, $i=1, \ldots, n$ correspond to the latent class
labels defined in (\ref{modexp}).

Similar to many adaptive methods, Lepki-type procedures suffer from a
dependency on a tuning parameter. In particular, in ERC method, a
constant governs the variance threshold (see \cite{Katkovnik99} or
\cite{ChichignoudLoustau13}), and in our selection rule as well, a
constant $b_1'>0$ appears in the majorant function of Lemma \ref{lemmamajkmeans}. As discussed earlier, the choice of this constant remains
an hard issue for application. In the sequel, we illustrate the
behavior of both adaptive methods w.r.t. 3 constants: $0.1$, $1$ and $10$.

\begin{figure}
\begin{tabular}{@{}cc@{}}

\includegraphics{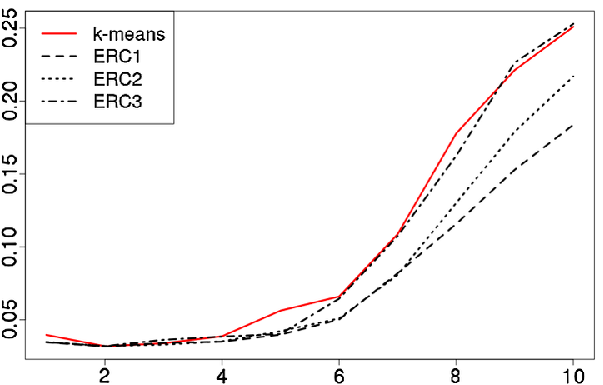}
 & \includegraphics{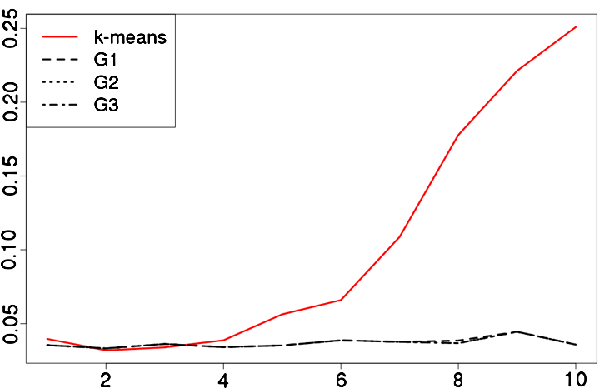}\\
\footnotesize{(a) $k$-means vs ERC method} & \footnotesize{(b) $k$-means vs Gradient}
\end{tabular}
\caption{Clustering risk averaged over 100 replications with $n=200$
for $k$-means versus ERC \textup{(a)} and the gradient \textup{(b)}.}\vspace*{-1pt}\label{fig1}
\end{figure}

Figure~\ref{fig1}(a)--(b) illustrates the evolution of the clustering risk
(\ref{cerror}) when $u\in\{1, \ldots, 10\}$ in model (\ref{modexp})
for $k$-means (red curve) versus both adaptive procedures.

In Figure~\ref{fig1}(a), we compare the clustering risk (\ref{cerror}) of
$k$-means (red curve) with ERC with 3 different constants (ERC1, ERC2
and ERC3). The methods are comparable, and we observe that ERC
performance is sensitive to the choice of the constant. Nevertheless, a
good calibration of this constant gives slightly better results than
$k$-means. In Figure~\ref{fig1}(b), the gradient approach with three different
constants (G1, G2 and G3) gives a clustering risk less than 5$\%$ for
any $u\in\{1, \ldots, 10\}$. In comparison, standard $k$-means
completely fails when $u$ is increasing. As a conclusion, our selection
rule significantly outperforms $k$-means and ERC for any constant. This
highlights the importance in practice to choose two different
bandwidths in each direction in this model, that is, an anisotropic
bandwidth. Our selection rule is also robust to the choice of the
constant, which confirms the theoretical study.

\begin{appendix}
\section*{Appendix}\label{appendix}

\setcounter{equation}{0}

\subsection{Proof of Lemma \texorpdfstring{\protect\ref{lemmadmargin}}{1}}
The proof is based on standard tools from differential calculus applied
to the
multivariate risk function $R\in\cC^2(U)$, where $U$ is an open ball
centered at $\f^\star$.
The first step is to apply a Taylor expansion of first order which
gives, for all $\f\in U$,
\begin{eqnarray*}
&& \R(\f)-\R\bigl(\f^\star\bigr)
\\
&&\qquad =\bigl(\f-\f^\star
\bigr)^\top\nabla\R\bigl(\f^\star \bigr)
\\
&&\quad\qquad{} +\sum
_{k\in\bN^m\dvtx \llvert  k\rrvert  =2} \frac{2(\f-\f^\star)^{k}}{k_1!\cdots
k_m!}\int_0^1
(1-t) \frac{\partial^2}{\partial\f^k}R\bigl(\f^\star +t\bigl(\f-\f^\star
\bigr)\bigr)\,dt,
\end{eqnarray*}
where
$\frac{\partial^2}{\partial\f^k}R=\frac{\partial^2}{\partial\f
_1^{k_1}\cdots\partial\f_m^{k_m}}
R$, $\llvert  k\rrvert  =k_1+\cdots+ k_m$ and $
(\f-\f^\star)^{k}=\prod_{j=1}^m(\f_j-\f^\star_j)^{k_j} $.
Now, by the property $ \nabla\R(\f^\star)=0 $ and the boundedness
of the second partial
derivatives,
we can write
\begin{eqnarray*}
\R(\f)-\R\bigl(\f^\star\bigr)&\leq&\kappa_1\sum
_{k\in\bN^m\dvtx \llvert  k\rrvert  =2} \bigl\llvert \f-\f ^\star\bigr\rrvert
^{k}\leq \kappa_1\sum_{i,j=1}^m
\bigl\llvert \f_i-\f^\star_i\bigr\rrvert \times
\bigl\llvert \f_j-\f^\star _j\bigr\rrvert
\\
&\leq& m\kappa_1 \bigl\llvert \f-\f^\star\bigr\rrvert
_2^2.
\end{eqnarray*}
It then remains to show the inequality
%
\begin{equation}
\label{secondstep} \bigl\llvert \f-\f^\star\bigr\rrvert _2\leq
2\bigl\llvert \D\bigl(\f,\f^\star\bigr)\bigr\rrvert _2/
\lambda_{\min},
\end{equation}
where $\lambda_{\min}$ is defined in the lemma. This can be done by using
standard inverse function theorem and the mean value theorem for
multi-dimensional functions.
Indeed, since the Hessian matrix of $ \R$---also viewed as the Jacobian
matrix of $ \D$---is positive definite at $\f^\star$, and since $\R
\in\cC^2(U)$, the inverse function
theorem shows the existence of a bijective function
$\D^{-1}\in\cC^1(G(U))$ such that
\[
\bigl\llvert \f-\f^\star\bigr\rrvert _2=\bigl\llvert
\D^{-1}\circ\D(\f)-\D^{-1}\circ\D\bigl(\f ^\star
\bigr)\bigr\rrvert _2\qquad\mbox{for any } \f\in U. %
\]
We can then apply a vector-valued version of the mean value
theorem to obtain
\begin{eqnarray}\label{mvt2}
\bigl\llvert \f-\f^\star\bigr\rrvert _2&\leq&
\sup_{u\in[G(\f),G(\f^\star)]}\duvvvert J_{\D
^{-1}}(u)\duvvvert _2
\bigl\llvert \D\bigl(\f^\star\bigr)-\D(\f)\bigr\rrvert _2
\nonumber\\[-8pt]\\[-8pt]
\eqntext{\mbox{for any } \f\in U,}
\end{eqnarray}
where $ [G(\f),G(\f^\star)] $ denotes the multi-dimensional bracket between
$G(\f)$ and $G(\f^\star)$, and $\vvvert \cdot\vvvert _2 $
denotes the operator norm associated to the Euclidean norm $\llvert  \cdot\rrvert  _2$.
Since $\llvert  \f-\f^\star\rrvert  _2\leq\delta$ and $G$ is continuous, we now have
\[
\lim_{\delta\to0}\sup_{u\in
[G(\f),G(\f^\star)]}\duvvvert
J_{\D^{-1}}(u)\duvvvert _2= \duvvvert J_{\D^{-1}} \bigl(G
\bigl(\f^\star\bigr)\bigr)\duvvvert _2. %
\]
Then, for $\delta>0$ small enough, we have with (\ref{mvt2})
\begin{eqnarray*}
\bigl\llvert \f-\f^\star\bigr\rrvert _2&\leq& 2\duvvvert
J_{\D^{-1}} \bigl(G\bigl(\f^\star\bigr)\bigr)\duvvvert
_2 \bigl\llvert \D\bigl(\f^\star\bigr)-\D(\f)\bigr\rrvert
_2
\\
&=& 2\duvvvert J^{-1}_{\D}\bigl(\f^\star\bigr)\duvvvert
_2 \bigl\llvert \D\bigl(\f^\star\bigr)-\D(\f)\bigr\rrvert
_2
\\
&=& 2\duvvvert H_{\R}^{-1}\bigl(\f^\star\bigr)\duvvvert
_2 \bigl\llvert \D\bigl(\f^\star\bigr)-\D(\f)\bigr\rrvert
_2,
\end{eqnarray*}
where $H_{\R}$ is the Hessian matrix of $R$. (\ref{secondstep})
follows easily, and the
proof is complete.

\subsection{Proofs of Section~\texorpdfstring{\protect\ref{sectionkmeans}}{3}}

\mbox{}

\begin{pf*}{Proof of Lemma \protect\ref{dclustering}}
The Hessian matrix of $\RC(\cdot)$ involves integrals over faces of
the Vorono\"i diagram $(V_j(\bc))_{j=1}^k$. For $i\neq j$, let us\vspace*{1pt}
denote the face (possibly empty) common to $V_i(\bc)$ and $V_j(\bc)$
as $F_{ij}$. Moreover, denote $\sigma(\cdot)$ the $(d-1)$-di\-mensional
Lebesgue measure. Then, since $f$ is continuous and $X\in[0,1]^d$,
uniform continuity arguments ensure that the integral $\int_{F_{ij}}\llvert  x-m\rrvert  _2^2f(x)\sigma(dx)$ exists and depends\vspace*{1pt} continuously on
the location of the center $m$, for any $i,j$ and for any $m\in\bR
^d$. Then we can use the following lemma due to \cite{pollard82}.

\begin{lemma}[(\cite{pollard82})]
\label{pollard}
Suppose $\mathbb{E}_P\llvert  X\rrvert  _2<\infty$ and $P$ has a continuous density
$f$ w.r.t. Lebesgue measure. Assume\vspace*{-1pt} integral $\int_{F_{ij}}\llvert  x-m\rrvert  _2^2f(x)\sigma(dx)$ exists and depends continuously on
the location of the centers, for any $i,j$ and for any $m\in\bR^d$.
Then if centers $c_i$, $i=1, \ldots, d$ are all distinct, $\RC(\cdot
)$ has a Hessian matrix $H_{\RC}(\cdot)$ made up of $d\times d$ blocks,
\begin{eqnarray*}
&& H_{\RC}(\bc) (i,j)
\\
&&\qquad =
\cases{\displaystyle 2\mathbb{P}\bigl(X\in V_i(\bc)
\bigr)-2\sum_{u\neq i}\delta_{iu}^{-1}
\int_{F_{iu}}f(x)\llvert x-c_i\rrvert
_2^2\sigma(dx), &\quad if $i=j$,
\vspace*{3pt}\cr
\displaystyle -2 \delta_{ij}^{-1}\int_{F_{ij}}f(x)
(x-c_i) (x-c_j)^\top\sigma (dx), &\quad
otherwise,}
\end{eqnarray*}
where $\delta_{ij}=\llvert  c_i-c_j\rrvert  _2$ and $\bc\in\bR^{dk}$.
\end{lemma}

Hence there exists $\delta>0$ such that $\RC(\cdot)\in\mathcal
{C}^2(U)$, and Lemma \ref{lemmadmargin} with $R=\RC$ completes the proof.
\end{pf*}

\begin{pf*}{Proof of Lemma \protect\ref{lemmamajkmeans}}
We start with the study of $\llvert  \DCn_h-\bE\DCn_h\rrvert  _{2,\infty}$. For
ease of
exposition, we denote by $P_n^Z$
the empirical measure with respect to $Z_i$, $i=1, \ldots, n$ and by
$P^Z$ the expectation w.r.t. the law of $Z$. Then we have
%
\begin{eqnarray}
\label{eqlemma1}
&& \llvert \DCn_h-\bE\DCn_h\rrvert
_{2,\infty}\nonumber
\\
&&\qquad =\sup_{\bc\in[0,1]^{dk}}\bigl\llvert \DCn
_h(\bc)-\bE\DCn_h(\bc)\bigr\rrvert _2
\\
&&\qquad \leq \sqrt{kd}\sup_{\bc,i,j}\biggl\llvert \bigl(P_n^Z-P^Z
\bigr) \biggl(\int_{V_j}2\bigl(x^i-c_j^i
\bigr)\widetilde K_h(Z-x)\,dx \biggr)\biggr\rrvert.\nonumber
\end{eqnarray}
The cornerstone of the proof is to apply a concentration inequality to this
supremum of empirical process. We use in the
sequel the following Talagrand-type
inequality; see, for example, \cite{ComteLacour13}.

\begin{lemma}
\label{lemmatal}
Let $\cX_1,\ldots, \cX_n$ be i.i.d. random variables, and let $\cS$
be a countable subset of
$\bR^{m}$.
Consider the random variable
\[
U_n(\cS):=\sup_{\bc\in\cS}\Biggl\llvert  \frac{1}{n}\sum_{l=1}^n\psi
_{\bc}(\cX_l)-\bE\psi_{\bc}
(\cX_l)\Biggr\rrvert   ,
\]
where $\psi_{\bc} $ is such that
$
\sup_{\bc\in\cS}\llvert  \psi_{\bc}\rrvert  _{\infty}\leq M$, $\bE U_n(\cS)\leq
E$ and $\sup_{\bc\in\cS}\bE [\psi_{\bc}(Z)^2 ]\leq v$.
Then, for any $\delta>0$, we have
\[
\bP \bigl(U_n(\cS)\geq (1+2\delta)E \bigr)\leq\exp \biggl(-
\frac{\delta^2nE}{6v} \biggr)\vee\exp \biggl(-\frac{
(\delta\wedge1)\delta nE}{21M} \biggr). %
\]
\end{lemma}

The proof of Lemma \ref{lemmatal} is omitted; see \cite
{ComteLacour13}. We hence have to compile the quantities $E,v$ and $M$
associated with the random
variable
\[
\widetilde \zeta_n=\sup_{\bc,i,j}\biggl\llvert
\bigl(P_n^Z-P^Z \bigr) \biggl(\int
_{V_j}2\bigl(x^i-c_j^i
\bigr)\widetilde K_h(Z-x)\,dx \biggr)\biggr\rrvert. %
\]
The compilation of $E:=E(h)>0$ uses the same path as \cite{ChichignoudLoustau13}, Lemma~3. More precisely, we can apply a chaining
argument to the function $\int_{V_j}2(x^i-u)\widetilde K_h(Z-x)\,dx$,
for any $u\in(0,1)$. Then we have, together with a maximum inequality
due to \cite{Massart07}, Chapter~6,
%
\begin{eqnarray}
\label{Hbound} \bE\widetilde\zeta_n\leq \frac{b_3}{2\sqrt{n}\Pi_h( \beta)}+
\frac{b_4}{2\sqrt{n}\Pi_h(
\beta+1/2)}\leq\frac{b_5}{\sqrt{n}
\Pi_h( \beta)}:=E(h),
\end{eqnarray}
where $\Pi_h( \beta):=\prod_{i=1}^d h_i^{ \beta_i}$ for $ \beta\in
\bR^d_+$ provided
that $\prod_{i=1}^dh_i^{-1/2}\geq b_1/b_1'$ (thanks to the definition of
$ \cH_a $ and $n$
sufficiently large). The constant $b_3,b_4,b_5>0$ can be explicitly
computed. This calculation is omitted for simplicity.
Besides, using \cite{ChichignoudLoustau13}, Lemma 1, with
$\psi_{\bc,i,j}(Z):=\int_{V_j}2(x^i-c_j^i)\widetilde K_h(Z-x)\,dx$, we have
%
\begin{equation}
\label{vbound} \sup_{\bc,i,j}\bE \bigl[\psi_{\bc,i,j}(Z)^2
\bigr]\leq\frac
{b_6}{\Pi_h(2 \beta)}:=v(h),
\end{equation}
whereas \cite{ChichignoudLoustau13}, Lemma 2, allows us to write
%
\begin{equation}
\label{Mbound} \sup_{\bc,i,j}\llvert \psi_{\bc,i,j}\rrvert
_{\infty}\leq\frac{b_7}{\Pi_h(
\beta+1/2)}:=M(h),
\end{equation}
where $b_6,b_7$ are absolute constants. Hence, Lemma \ref{lemmatal},
together with (\ref{eqlemma1})--(\ref{Mbound}), gives us, for all $
\delta>0$,
\begin{eqnarray*}
&& \bP \bigl(\llvert \DCn_h-\bE\DCn_h\rrvert
_{2,\infty}\geq \sqrt{kd}(1+2\delta)E(h) \bigr)
\\
&&\qquad \leq\exp \biggl(-
\frac{\delta
^2nE(h)}{6v(h)} \biggr)\vee\exp \biggl(-\frac{(\delta\wedge1)\delta
nE(h)}{21M(h)} \biggr). %
\end{eqnarray*}
Moreover, note that from the previous calculations, we have
$nE(h)/v(h)=c\sqrt{n}/\Pi_h(
\beta)$
and $nE(h)/M(h)=c'\sqrt{n}\sqrt{\Pi_h(1/2)}$, where $c,c'>0$ depend
on $b_5,b_6$ and $b_5,b_7$, respectively. Provided that $\sqrt{n}(c\Pi
_h( \beta)\wedge
c'\sqrt{\Pi_h(1/2)})\geq(\log n)^2$ (thanks to the definition of $
\cH_a $ and $n$
sufficiently large), we come up with
\begin{eqnarray*}
&& \bP \bigl(\llvert \DCn_h-\bE\DCn_h\rrvert
_{2,\infty}\geq \sqrt{kd}(1+2\delta)E(h) \bigr)
\\
&&\qquad \leq\exp \biggl\{- \biggl(
\frac
{\delta^2}{6}\wedge\frac{
(\delta\wedge1)\delta}{21} \biggr) (\log n)^2 \biggr
\}.
\end{eqnarray*}
This gives us the first part of the majorant of Lemma \ref{lemmamajkmeans}.

The last\vspace*{1pt} step is to show a similar bound for the auxiliary empirical process
$\llvert  \DCn_{h,\eta}-\bE\DCn_{h,\eta}\rrvert  _{2,\infty}$. This can be easily
done by using Lemma
\ref{lemmatal} together with the previous results. Then we have for
any $h,\eta\in\cH_a$,
\begin{eqnarray*}
&& \bP \bigl(\llvert \DCn_{h,\eta}-\bE\DCn_{h,\eta}\rrvert
_{2,\infty}\geq \sqrt{kd}(1+2\delta)E(h\vee\eta) \bigr)
\\
&&\qquad \leq\exp \biggl\{-
\biggl(\frac{\delta^2}{6}\wedge\frac{
(\delta\wedge1)\delta}{21} \biggr) (\log n)^2
\biggr\}, %
\end{eqnarray*}
where with a slight abuse of notation, the maximum $\vee$ is
understood coordinatewise.
Using the union bound, the definition of $\mathcal{M}^{\mathrm
{k}}_l(\cdot,\cdot)$ allows us to
write
\begin{eqnarray*}
&& \bP \Bigl(\sup_{h,\eta} \bigl\{\llvert \DCn_{h,\eta}-\bE
\DCn_{h,\eta
}\rrvert _{2,\infty}+\llvert \DCn_h-\bE
\DCn_h\rrvert _{
2,\infty} -\mathcal{M}^{\mathrm{k}}_l(h,
\eta) \bigr\}> 0 \Bigr)
\\
&&\qquad \leq (\operatorname{card}\cH_a )^2 \sup_{h,\eta}
\bP \bigl(\llvert \DCn_{h,\eta}-\bE\DCn_{h,\eta
}\rrvert
_{2,\infty}
\\
&&\hspace*{107pt}{}+\llvert \DCn_h-\bE\DCn_h\rrvert
_{2,\infty} -\mathcal{M}^{\mathrm{k}}_l(h,\eta)> 0 \bigr)
\\
&&\qquad \leq (\operatorname{card}\cH_a )^2 \sup_{h,\eta}
\bigl\{\bP \bigl(\llvert \DCn_h-\bE\DCn_h\rrvert
_{2,\infty}- \sqrt{kd}(1+2\delta)E(h)>0 \bigr)
\\
&&\hspace*{98pt}{}+  \bP \bigl(\llvert
\DCn_{h,\eta}-\bE \DCn_{h,\eta}\rrvert _{2,\infty}
\\
&&\hspace*{185pt}{}-
\sqrt{kd}(1+2\delta)E(h\vee\eta)>0 \bigr) \bigr\}
\\
&&\qquad \leq2 (\operatorname{card}\cH_a )^2\exp \biggl(-
\frac{\delta
^2}{6}\wedge\frac{(\delta\wedge
1)\delta}{21}(\log n)^2 \biggr)\leq
n^{-l},
\end{eqnarray*}
where we choose $b'_1=b_5(1+2\delta)$ with $\delta:=\delta(l)=1\vee
(21(l+2)/(\log n))$.
\end{pf*}

\begin{pf*}{Proof of Theorem \protect\ref{thmkmeans}}
The proof of Theorem \ref{thmkmeans} is a direct application of
Theorem \ref{thmainresult} and Lemma \ref{lemmamajkmeans}. Indeed, for any $l\in\bN
^\star$, for $n$ large enough, we have with probability
$1-n^{-l}$,
\[
\bigl\llvert \DC\bigl(\bcn_{\widehat h},\bc^\star\bigr)\bigr\rrvert
_2\leq 3\inf_{h\in\cH_a} \bigl\{\bias(h)+
\cM^{\mathrm{k},\infty
}_l(h) \bigr\}, %
\]
where $\bias(h)$ is defined as
\[
\bias(h):=\max \Bigl(\llvert \bE\DCn_{h}-\DC\rrvert _{2,\infty},
\sup_{\eta
}\llvert \bE\DCn_{h,\eta}-\bE
\DCn_\eta\rrvert _{2,
\infty} \Bigr)\qquad \forall h\in
\cH_a. %
\]
The control of the bias function is as follows:
\begin{eqnarray*}
&& \llvert \bE\DCn_{h,\eta}-\bE\DCn_\eta\rrvert
_{2,\infty}^2
\\
&&\qquad = \sup
_{\bc\in[0,1]^{dk}}\sum_{i,j} \biggl\{ \int
_{V_j}2\bigl(x^i-c_j^i
\bigr) \bigl(\bE_{P^Z} \widetilde K_{h,\eta}(Z-x)-
\bE_{P^Z}\widetilde K_\eta(Z-x) \bigr)\,dx \biggr
\}^2
\\
&& \qquad =\sup_{\bc\in[0,1]^{dk}}\sum
_{i,j} \biggl\{\int_{V_j}2
\bigl(x^i-c_j^i\bigr) \bigl(
\bE_{P^X} K_{h,\eta}(X-x)-\bE_{P^X}
K_\eta(X-x) \bigr)\,dx \biggr\}^2
\\
&&\qquad \leq 4\sup_{\bc\in[0,1]^{dk}}\sum_{i,j}\int
_{V_j}\bigl(x^i-c_j^i
\bigr)^2\,dx\bigl\llvert K_{\eta}*(K_h*f-f)\bigr
\rrvert _2^2
\\
&&\qquad \leq 4k\bigl\llvert \mathcal{F}[K]\bigr\rrvert _{\infty}\llvert
f_h-f\rrvert _2^2,
\end{eqnarray*}
where $\llvert  f_h-f\rrvert  _2:=\llvert  K_h*f-f\rrvert  _2$ is the usual nonparametric bias term
in deconvolution estimation. Besides, note that
\begin{eqnarray*}
&& \llvert \bE\DCn_{h}-\DC\rrvert _{2,\infty}^2
\\
&&\qquad =\sup
_{\bc\in[0,1]^{dk}}\sum_{i,j} \biggl\{\int
_{V_j}2\bigl(x^i-c_j^i
\bigr) \bigl(\bE_{P^X} K_{h}(X-x)-f(x) \bigr)\,dx \biggr
\}^2
\\
&&\qquad \leq 4\sup_{\bc\in[0,1]^{dk}}\sum_{i,j}\int
_{V_j}\bigl(x^i-c_j^i
\bigr)^2\,dx\llvert K_{h}*f-f\rrvert _2^2.
\end{eqnarray*}
Then we need a control of the bias function,
\[
B^{\mathrm{k}}(h):=2\sqrt{k} \bigl(1\vee\bigl\llvert \cF[K]\bigr\rrvert
_{\infty} \bigr)\llvert K_h*f-f\rrvert _2\qquad \forall h\in\cH.
\]
By using Comte and Lacour \cite{ComteLacour13}, Proposition~3, we directly have for all
$f\in\mathcal{N}_{2,d}(s,L)$,
%
\begin{equation}
\label{bcontrolkmeans} B^{\mathrm{k}}(h)\leq2\sqrt{k} \bigl(1\vee\bigl\llvert \cF[K]
\bigr\rrvert _{\infty} \bigr)L\sum_{j=1}^dh_j^{s_j}\qquad \forall h\in\cH.
\end{equation}
Now, we have to use a result such as Lemma \ref{dclustering}, for our
family of estimators $\{\bcn_{h}, h\in\mathcal{H}_a\}$. In other
words, we need to check that this family of estimators is consistent
with respect to the Euclidean norm in $\mathbb{R}^{dk}$.

\begin{lemma}
\label{consistencykmeans}
Assume $f$ is continuous, $X\in[0,1]^d$ a.s. and the Hessian matrix of
$\RC$ is positive definite on $\cM$. Consider the family $\{\bcn
_{h},h\in\cH_a\}$ with $\cH_a$ defined in Lemma \ref{lemmamajkmeans}. Then, for any $\delta>0$, for any $l\in\mathbb{N}^\star
$, for any $\bcn_h\in\cH_a$, there exists $\bc^\star\in\mathcal
{M}$ such that for $n$ great enough, with probability $1-n^{-l}$,
\begin{eqnarray*}
\bigl\llvert \bcn_h-\bc^\star\bigr\rrvert _2
\leq\delta.
\end{eqnarray*}
\end{lemma}

\begin{pf} 
Using \cite{gg},
the positive definiteness of the Hessian matrix on $\mathcal{M}$ and
the continuity of $f$, we have, for any $\bcn_h\in\cH_a$, for some
constant $A_1>0$, $\llvert  \bcn_h-\bc^\star\rrvert  _2\leq A_1(\RC(\bcn_h)-\RC
(\bc^\star))$, where $\bc^\star\in\arg\min_{\bc\in\cM}\llvert  \bcn
_h-\bc\rrvert  _2$. It remains to show that by definition of $\cH_a$ in Lemma
\ref{lemmamajkmeans}, with high probability, $\RC(\bcn_h)-\RC(\bc
^\star)\to0$ as $n$ tends to infinity. This can be seen easily from
Chichignoud and Loustau \cite{ChichignoudLoustau13}, which gives the
order of the bias term and the variance term for such a problem. At
this stage, we can notice that localization is used in \cite
{ChichignoudLoustau13}, and appears to be necessary here. However,
using a global approach (i.e., a simple Hoeffding inequality to the
family of kernel ERM), we can have, for any $l\in\mathbb{N}^\star$,
with probability $1-n^{-l}$,
\[
\RC(\widehat\bc_h)-\RC\bigl(\bc^\star\bigr)\lesssim
\frac{\Pi_h(-\beta
)}{\sqrt{n}}+\sum_{j=1}^dh_j^{s_j}\qquad \forall h\in\cH_a. %
\]
By definition of $\cH_a$, the RHS tends to zero as $n\to\infty$, and
then for $n$ great enough, this term is controlled by $\delta$. 
\end{pf}

Then, for any $h\in\cH_a$ and $n$ great enough, Lemma \ref
{dclustering} allows us to write with probability $1-n^{-l}$,
\[
\sqrt{\RC(\bcn_h)-\RC\bigl(\bc^\star\bigr)}\leq 2
\frac{\sqrt{kd}}{\lambda_{\min}}\bigl\llvert \nabla\RC\bigl(\bcn_h,
\bc^\star\bigr)\bigr\rrvert _2.
\]
Using Theorem \ref{thmainresult} with $ l=q $, bias control (\ref
{bcontrolkmeans}) and the last inequality, there
exists an absolute constant $ b_8>0 $ such that
\[
\sup_{f\in\mathcal{N}_{2}(s,L)}\bE \bigl[\RC(\bcn_{\widehat
h})-\RC\bigl(
\bc^\star\bigr) \bigr]\leq b_8 \inf_{h\in\cH_a}
\Biggl\{\sum_{j=1}^dh_j^{ s_j}+
\frac{\Pi_h(-\beta
)}{n} \Biggr\} ^2+b_8n^{-q}.
\]
Let $ h^\star$ denote the oracle bandwidth as
$
h^\star:=\arg\inf_{h\in\cH} \{\sum_{j=1}^dh_j^{ s_j}+\frac
{\Pi_h(-\beta)}{n} \}$,
and define the oracle bandwidth $ h^\star_a $ on the net $\cH_a$ such
that $
a h^\star_{a,j}\leq h^\star_{j}\leq h^\star_{a,j} $, for all $
j=1,\dots,d$.
Eventually, we have
\[
\sup_{f\in\mathcal{N}_{2}(s,L)}\bE \bigl[\RC(\bcn_{\widehat
h})-\RC\bigl(
\bc^\star\bigr) \bigr]\leq b_8 a^{-qd/2}\inf
_{h\in\cH} \Biggl\{\sum_{j=1}^dh_j^{ s_j}+
\frac{\Pi
_h(-\beta)}{n} \Biggr\}^2+b_8n^{-q}.
\]
By a standard bias variance trade-off, we obtain the assertion of the
theorem, provided that $q\geq1$.
\end{pf*}

\subsection{Proofs of Section~\texorpdfstring{\protect\ref{sectionlocalglobal}}{4}}

\mbox{}

\begin{pf*}{Proof of Lemma \protect\ref{lemlocalmargincondition}}
By definition, we first note that
\[
\bigl\llvert {G^{\mathrm{loc}} \bigl(\estim_{h}(x_0)
\bigr) -G^{\mathrm{loc}} \bigl(f^\star(x_0) \bigr)}\bigr
\rrvert =\bigl\llvert \bE\rho_\gamma ' \bigl(\xi
_1+f^\star(x_0)-\estim_{
h}(x_0)
\bigr)-\bE\rho_\gamma' (\xi_1 )\bigr\rrvert.
\]
Using the mean value theorem and the assumption
$
\sup_{h\in\cH}\llvert  \estim_{h}(x_0)-f^\star(x_0)\rrvert  \leq\bE\rho_\gamma
''(\xi)/4
$, there exists $ c \in[-\bE\rho_\gamma''(\xi_1)/4,\bE\rho
_\gamma''(\xi_1)/4]$ such that
\[
\bigl\llvert {G^{\mathrm{loc}} \bigl(\estim_{h}(x_0)
\bigr) -G^{\mathrm{loc}} \bigl(f^\star(x_0) \bigr)}\bigr
\rrvert ={\bE\rho_\gamma ''(\xi
_1+c)}\bigl\llvert f^\star(x_0)-
\estim_{h}(x_0)\bigr\rrvert. %
\]
Since $ \bE\rho_\gamma''(\xi_1+\cdot) $ is a 2-Lipschitz function,
it yields
\[
\bigl\llvert {G^{\mathrm{loc}} \bigl(\estim_{h}(x_0)
\bigr) -G^{\mathrm{loc}} \bigl(f^\star(x_0) \bigr)}\bigr
\rrvert \geq \frac{\bE\rho_\gamma''(\xi_1)}{2}\bigl\llvert f^\star(x_0)-
\estim_{h}(x_0)\bigr\rrvert. %
\]
The proof is complete.
\end{pf*}

%
%

\begin{pf*}{Proof of Theorem \protect\ref{thholderadapation}}
From \cite{ChichignoudLederer13}, Theorem 1, we notice that all
estimators $ \{\estim_h(x_0),h\in\cH\} $ are consistent, and thus,
for $n$ sufficiently large, the assumption
of Lemma \ref{lemlocalmargincondition} holds for all $ x_0\in\cT$.
Using Theorem \ref{thmainresult} with $ l>0 $ and Lemma \ref{lemlocalmargincondition}, we get
\[
\bigl\llvert \estim_{\widehat
h^\mathrm{loc}}(x_0)-f^\star(x_0)
\bigr\rrvert \leq\frac{6}{\bE\rho_\gamma
''(\xi_1)}\inf_{
h\in\cH_a} \bigl\{
\bias(h)+2\cM_{l}^{\mathrm{loc},\infty}(h) \bigr\},
\]
with
$
\bias(h)=\max (\llvert  \bE\widehat{G}^{\mathrm{loc}}_{h}-G^{\mathrm
{loc}}\rrvert  _{\infty},\sup_{\eta\in
\cH}\llvert  \bE\widehat{G}^{\mathrm{loc}}_{h,\eta}
-\bE\widehat{G}^{\mathrm{loc}}_\eta\rrvert  _
{
\infty} )$.
The control of $\bias(\cdot)$ over H\"{o}lder classes is based on the
same schema
as in
\cite{GoldenshlugerLepski08}, applied to the function
$
F_t(\cdot):=\bE\rho_\gamma'(f^\star(\cdot)-t+\xi_1)$.
For any $
f\in\Sigma({ s},L) $ and any $ h\in\cH$, we then want
to show
%
\begin{eqnarray}\label{eqlocalbiascontrol}
B^{\mathrm{loc}}(h)&\leq&\sup_{t\in[-B,B]}\sup
_{y\in\cT}\biggl\llvert \int K_{h}(x-y)
\bigl[F_t(x)-F_t(y) \bigr]\,dx\biggr\rrvert
\nonumber\\[-8pt]\\[-8pt]\nonumber
&\leq& L
\llvert K\rrvert _\infty\sum_{j=1}^dh_j^{ s_j}.
\end{eqnarray}
By definition, we see that
$
\llvert  \bE\widehat{G}^{\mathrm{loc}}_{h}-G^{\mathrm{loc}}\rrvert  _{\infty}=\sup_{t\in[-B,B]}\llvert  \bE
K_{h}(W-x_0) [F_t(W)-F_t(x_0) ]\rrvert
$
and by definition of $ \bE\widehat{G}^{\mathrm{loc}}_{h,\eta} $ and
$F_t$, we have
\begin{eqnarray*}
-\bE\widehat{G}^{\mathrm{loc}}_{h,\eta}(t)&=&\int F_t(x)K_{h,\eta}(x-x_0)\,dx
\\
&=&\int F_t(x) \biggl(\int K_{h}(x-y)K_{\eta}(y-x_0)\,dy
\biggr)\,dx.
\end{eqnarray*}
Using Fubini's theorem and the equation $ \int K_{h}(x-y)\,dx=1 $ for all
$ y\in\cT$, we get
\begin{eqnarray*}
-\bE\widehat{G}^{\mathrm{loc}}_{h,\eta}(t)&=&\int K_{\eta}(y-x_0)
F_t(y)\,dy
\\
&&{} +\int K_{\eta
}(y-x_0) \biggl( \int
K_{h}(x-y) \bigl[F_t(x)-F_t(y) \bigr]\,dx
\biggr)\,dy
\\
&=&\int K_{\eta}(y-x_0) F_t(y)\,dy
\\
&&{}+\int
K_{\eta}(y-x_0) \int K_{h}(x-y)
\bigl[F_t(x)-F_t(y) \bigr]\,dx\,dy.
\end{eqnarray*}
Then it holds for any $ x_0\in\cT$,
\begin{eqnarray*}
&& \bigl\llvert \bE\widehat{G}^{\mathrm{loc}}_{h,\eta}(t)-\bE\widehat
{G}^{\mathrm{loc}}_{\eta}(t)\bigr\rrvert
\\
&&\qquad =\biggl\llvert \int
K_{\eta}(y-x_0) \int K_{h}(x-y)
\bigl[F_t(x)-F_t(y) \bigr]\,dx\,dy\biggr\rrvert
\\
&&\qquad \leq \bigl\llVert K_{\eta}(\cdot-x_0)\bigr\rrVert
_1\sup_{y\in\cT}\biggl\llvert \int
K_{h}(x-y) \bigl[F_t(x)-F_t(y) \bigr]\,dx
\biggr\rrvert
\\
&&\qquad =\sup_{y\in\cT}\biggl\llvert \int K_{h}(x-y)
\bigl[F_t(x)-F_t(y) \bigr]\,dx\biggr\rrvert.
\end{eqnarray*}
We have then shown the first inequality in (\ref{eqlocalbiascontrol}). Using the smoothness
of $ \rho_\gamma' $, we have
for all $
f\in\Sigma( s,L)
$,
\begin{eqnarray*}
&& \biggl\llvert \int K_{h}(x-y) \bigl[F_t(x)-F_t(y)
\bigr]\,dx\biggr\rrvert
\\
&&\qquad =\biggl\llvert \int K_{h}(x-y)\bE \bigl[
\rho_\gamma' \bigl(f(x)-t+\xi_1 \bigr)-
\rho_\gamma' \bigl(f(y)-t+\xi _1 \bigr) \bigr]
\,dx\biggr\rrvert
\\
&&\qquad \leq\biggl\llvert \int K_h(x-y) \bigl(f(x)-f(y)\bigr)\,dx\biggr
\rrvert
\\
&&\qquad \leq L\llVert K\rrVert _\infty\sum_{j=1}^dh_j^{ s_j}.
\end{eqnarray*}
Therefore, (\ref{eqlocalbiascontrol}) holds. Then, using Theorem
\ref{thmainresult} with $ l=q $, Lemma \ref{lemlocalmargincondition} and (\ref{eqlocalbiascontrol}), there
exists an absolute constant $ T_1>0 $ such that
\[
\sup_{f\in\Sigma({ s},L)}\bE\bigl\llvert \estim_{\widehat
h}(x_0)-f(x_0)
\bigr\rrvert ^q\leq T_1 \inf_{h\in\cH_a}
\Biggl\{\sum_{j=1}^dh_j^{ s_j}+
\sqrt{\frac{\log
(n)}{n\Pi_h}} \Biggr\} ^q+T_1n^{-q}.
\]
Let $ h^\star$ denote the oracle bandwidth as
$
h^\star:=\arg\inf_{h\in\cH} \{\sum_{j=1}^dh_j^{ s_j}+\sqrt
{\frac{\log(n)}{n\Pi_h}} \}$,
and define the oracle bandwidth $ h^\star_a $ such that $
a h^\star_{a,j}\leq h^\star_{j}\leq h^\star_{a,j} $, for all $
j=1,\dots,d$.
Then we get
\[
\sup_{f\in\Sigma({ s},L)}\bE\bigl\llvert \estim_{\widehat
h}(x_0)-f(x_0)
\bigr\rrvert ^q\leq T_1 a^{-qd/2}\inf
_{h\in\cH} \Biggl\{\sum_{j=1}^dh_j^{ s_j}+
\sqrt{\frac
{\log(n)}{n\Pi_h}} \Biggr\} ^q+T_1n^ { -q }.
\]
By a standard bias variance trade-off, we obtain the assertion of the
theorem.
\end{pf*}

\begin{pf*}{Proof of Theorem \protect\ref{thnikolskiiadapation}}
Here again, the assumption of Lemma \ref{lemlocalmargincondition} holds for $n$ sufficiently large for all $ x_0\in\cT$.
Using Theorem \ref{thmainresult} with $ l>0 $ and adding the $ \bL
_q $-norm, we have
\[
\llVert \estim_{\widehat
h^{\mathrm{glo}}_q}-f\rrVert _q\leq\frac{6}{\bE\rho_\gamma''(\xi
_1)}\inf
_{h\in\cH} \bigl\{ B(h)+2\Gamma_{l,q}^{\mathrm{glo}}(h)
\bigr\}, %
\]
where
$
B(h)=\max (\llVert  \bE\widehat{G}^{\mathrm{loc}}_{h}-G^{\mathrm
{loc}}\rrVert  _{q,\infty},\sup_{\eta\in\cH
}\llVert  \bE\widehat{G}^{\mathrm{loc}}_{h,
\eta}
-\bE\widehat{G}^{\mathrm{loc}}_\eta\rrVert  _{q,
\infty} )$.
The control of the bias term is based on the schema of
\cite{GoldenshlugerLepski11}
for linear estimates. For any $ h\in\cH$, we
want
to show that
%
\begin{equation}
\label{eqglobalbiascontrol} B(h)\leq\sup_{t\in[-B,B]}\biggl\llVert \int
K_{h}(x-\cdot) \bigl[F_t(x)-F_t(\cdot)
\bigr]\,dx\biggr\rrVert _q \leq L\sum_{j=1}^dh_j^{ s_j},
\end{equation}
where we recall $ F_t(x):=\bE\rho_\gamma'(f(x)-f_t(x)+\xi_1)$.
By definition, one has
\[
\bigl\llVert \bE\widehat{G}^{\mathrm{loc}}_{h}-G^{\mathrm{loc}}\bigr
\rrVert _{q,\infty
}=\sup_{t\in[-B,B]}\bigl\llVert \bE
K_{h}(W-\cdot) \bigl[F_t(W)-F_t(\cdot)
\bigr]\bigr\rrVert _q. %
\]
%
Moreover, in the proof of
Theorem \ref{thholderadapation}, we
have shown that for any $ x_0\in\cT$,
\begin{eqnarray*}
&& \bE\widehat{G}^{\mathrm{loc}}_{\eta}(t,x_0)-\bE\widehat
{G}^{\mathrm{loc}}_{h,\eta}(t,x_0)
\\
&&\qquad =\int K_{\eta
}(y-x_0)
\int K_{h}(x-y) \bigl[F_t(x)-F_t(y)
\bigr]\,dx\,dy. %
\end{eqnarray*}
By Young's inequality and the definition of the kernel in Section~\ref{sectionKERM}, it yields
\begin{eqnarray*}
&& \bigl\llVert \bE\widehat{G}^{\mathrm{loc}}_{\eta}-\bE
\widehat{G}^{\mathrm
{loc}}_{h,\eta}\bigr\rrVert _{q,\infty}
\\
&&\qquad =\sup _{t\in
[-B,B]}\biggl\llVert \int K_{\eta}(y-\cdot) \int
K_{h}(x-y) \bigl[F_t(x)-F_t(y) \bigr]\,dx\,dy
\biggr\rrVert _{q,\infty}
\\
&&\qquad \leq\sup_{t\in[-B,B]}\biggl\llVert \int K_{h}(x-\cdot)
\bigl\llvert F_t(x)-F_t(\cdot)\bigr\rrvert\, dx\biggr
\rrVert _{q,\infty}.
\end{eqnarray*}
Using the smoothness of $ \rho_\gamma' $, we have for any $ x,y\in
\cT$ and any $
{t\in[-B,B]} $,
\begin{eqnarray*}
&& F_t(x)-F_t(y)=\bE \bigl[ \rho_\gamma'
\bigl(f(x)-t+\xi_1 \bigr)-\rho_\gamma'
\bigl(f(y)-t+\xi _1 \bigr) \bigr] \leq\bigl\llvert f(x)-f(y)\bigr
\rrvert.
\end{eqnarray*}
Therefore, (\ref{eqglobalbiascontrol}) holds for all $ f\in\cN
_{q,d}( s,L)
$. Then,
using Theorem \ref{thmainresult} with $ l=q $, Lemma \ref{lemlocalmargincondition} and (\ref{eqglobalbiascontrol}),
there
exists an absolute constant $ T_2>0 $ such that
\begin{eqnarray*}
&& \sup_{f\in\cN_{q,d}({ s},L)}\bE\llVert \estim_{\widehat
h^{\mathrm{glo}}_q}-f\rrVert
_q^q
\\
&&\qquad \leq T_2\times\cases{\displaystyle \inf
_{h\in\cH} \Biggl\{\sum_{j=1}^dh_j^{
s_j}+(n
\Pi_h)^{-(q-1)/q} \Biggr\}^q+n^{-q}, &\quad
if $q\in[1,2[$,
\cr
\displaystyle\inf_{h\in\cH} \Biggl\{\sum
_{j=1}^dh_j^{ s_j}+(n
\Pi_h)^{-1/2} \Biggr\}^q+n^{-q}, &\quad
if $q\in[2,\infty[$.}
\end{eqnarray*}
Computing these infimums, we obtain the assertion of the theorem.
\end{pf*}
\end{appendix}



%

\printaddresses
\end{document}